# Derivative of Map of Banach algebra

Aleks Kleyn

Aleks_Kleyn@MailAPS.org
http://AleksKleyn.dyndns-home.com:4080/
http://sites.google.com/site/AleksKleyn/
http://arxiv.org/a/kleyn_a_1
http://AleksKleyn.blogspot.com/

Abstract. Let $A$ be Banach algebra over commutative ring $D$. The map $f:A\rightarrow A$ is called differentiable in the Gâteaux sense, if

$$f(x+a)-f(x)=\partial f(x)\circ a+o(a)$$

where the Gâteaux derivative $\partial f(x)$ of map $f$ is linear map of increment $a$ and $o$ is such continuous map that

$$\lim_{a\rightarrow 0}\frac{|o(a)|}{|a|}=0$$

Assuming that we defined the Gâteaux derivative $\partial^{n-1}f(x)$ of order $n-1$, we define

$$\partial^n f(x)\circ(a_1\otimes...\otimes a_n)=\partial(\partial^{n-1}f(x)\circ(a_1\otimes...\otimes a_{n-1}))\circ a_n$$

the Gâteaux derivative of order $n$ of map $f$. Since the map $f(x)$ has all derivatives, then the map $f(x)$ has Taylor series expansion

$$f(x)=\sum_{n=0}^{\infty}(n!)^{-1}\partial^n f(x_0)\circ(x-x_0)^n$$

# Contents

CHAPTER 1

# Preface

## 1.1. Preface

The possibility of linear approximation of a map is at the heart of calculus and main constructions of calculus have their roots in linear algebra. Since the product in the field is commutative, then linear algebra over a field is relatively simple. When we explore algebra over commutative ring, we still see some familiar statements of linear algebra; however, we meet new statements, which change the landscape of linear algebra.

Here I want to draw attention to the evolution of the concept of the derivative from the time of Newton. When we study functions of single variable, the derivative in selected point is a number

$$dx^2 = 2x\ dx$$

When we study function of multiple variables, we realize that it is not enough to use number. The derivative becomes vector or gradient

$$z = x^2 + y^3$$

$$dz = 2x\ dx + 3y^2\ dy$$

When we study maps of vector spaces, this is a first time that we tell about derivative as operator

$$x = u\sin v \qquad y = u\cos v \qquad z = u$$

$$dx = \sin v\ du + u\cos v\ dv \quad dy = \cos v\ du - u\sin v\ dv \quad dz = 1\ du + 0\ dv$$

However since this operator is linear, then we can represent derivative as matrix. Again we express a vector of increment of a map as product of a matrix of derivative (Jacobian matrix) over vector of increment of argument

$$\begin{pmatrix}dx\\ dy\\ dz\end{pmatrix} = \begin{pmatrix}\sin v & u\cos v\\ \cos v & -u\sin v\\ 1 & 0\end{pmatrix}\begin{pmatrix}du\\ dv\end{pmatrix}$$

Since the algebra is a module over some commutative ring, there exist two ways to explore structures generated over the algebra.

If the algebra is a free module, then we can choose a basis and consider all operations in coordinates relative a given basis. Although the basis can be arbitrary, we can choose the simplest basis in terms of algebraic operations. Beyond doubt, this approach has the advantage that we are working in commutative ring where all operations are well studied.

Exploration of operations in algebra regardless of the chosen basis gives an opportunity to consider elements of algebra as independent objects. I considered the structure of linear map of algebra in the book [1]. In this book, I used this tool to study calculus over noncommutative algebra.

## 1.2. Conventions

CONVENTION 1.2.1. *I assume sum over index $i$ in expression like*

$$a_{i\cdot 0}xa_{i\cdot 1}$$

□

CONVENTION 1.2.2. *Since the tensor* $a\in A^{\otimes(n+1)}$ *has the expansion*

$$a=a_{i\cdot 0}\otimes a_{i\cdot 1}\otimes\dots\otimes a_{i\cdot n}\quad i\in I$$

*then set of permutations* $\sigma=\{\sigma_i\in S(n):i\in I\}$ *and tensor $a$ generate the map*

$$(a,\sigma):A^{\times n}\to A$$

*defined by rule*

$$\begin{aligned}(a,\sigma)\circ(b_1,...,b_n)&=(a_{i\cdot 0}\otimes a_{i\cdot 1}\otimes\dots\otimes a_{i\cdot n},\sigma_i)\circ(b_1,...,b_n)\\&=a_{i\cdot 0}\sigma_i(b_1)a_{i\cdot 1}...\sigma_i(b_n)a_{i\cdot n}\end{aligned}$$

□

CONVENTION 1.2.3. *Let the tensor* $a\in A^{\otimes(n+1)}$ . *When* $x_1=...=x_n=x$, *we assume*

$$a\circ x^n=a\circ(x_1\otimes\dots\otimes x_n)$$

□

CONVENTION 1.2.4. *Element of $\Omega$-algebra $A$ is called* $A$**-number**. *For instance, complex number is also called $C$-number, and quaternion is called $H$-number.* □

CONVENTION 1.2.5. *Let $A$ be associative $D$-algebra. The representation*

$$A\otimes A \xrightarrow[*]{f} A\qquad f(p):a\to p\circ a$$

*of $D$-module* $A\otimes A$ *is defined by the equation*

$$(a\otimes b)\circ c=acb \tag{1.2.1}$$

*and generates the set of linear maps. This representation generates product $\circ$ in $D$-module* $A\otimes A$ *according to rule*

$$(p\circ q)\circ a=p\circ(q\circ a)$$

□

CHAPTER 2

# Differentiable Maps

In this chapter, we explore derivative and differential of the map into $D$-algebra. Complex field is the algebra over real field. In the calculus of functions of complex variable, we consider linear maps generated by the map

$$I_0 \circ z = z$$

In this section, we also consider linear maps generated by the map

$$I_0 \circ a = a$$

## 2.1. Topological Ring

Definition 2.1.1. *Ring $D$ is called* **topological ring**[2.1] *if $D$ is topological space and the algebraic operations defined in $D$ are continuous in the topological space $D$.* □

According to definition, for arbitrary elements $a$, $b \in D$ and for arbitrary neighborhoods $W_{a-b}$ of the element $a-b$, $W_{ab}$ of the element $ab$ there exists neighborhoods $W_a$ of the element $a$ and $W_b$ of the element $b$ such that $W_a - W_b \subset W_{a-b}$, $W_a W_b \subset W_{ab}$.

Definition 2.1.2. **Norm on ring** *$D$ is a map*[2.2]

$$d \in D \rightarrow |d| \in R$$

*which satisfies the following axioms*

- $|a| \geq 0$
- $|a| = 0$ *if, and only if,* $a = 0$
- $|ab| = |a|\ |b|$
- $|a+b| \leq |a| + |b|$

*Ring $D$, endowed with the structure defined by a given norm on $D$, is called* **normed ring**. □

Invariant distance on additive group of ring $D$

$$d(a, b) = |a - b|$$

defines topology of metric space, compatible with ring structure of $D$.

Definition 2.1.3. *Let $D$ be normed ring. Element $a \in D$ is called* **limit of a sequence** $\{a_n\}$

$$a =$$

2.1 I made definition according to definition from [4], chapter 4

2.2 I made definition according to definition from [2], IX, §3.2 and definition [6]-1.1.12, p. 23.

*if for every $\epsilon \in R$, $\epsilon > 0$, there exists positive integer $n_0$ depending on $\epsilon$ and such, that*

$$|a_n - a| < \epsilon$$

*for every $n > n_0$.* □

Theorem 2.1.4. *Let $D$ be normed ring of characteristic 0 and let $d \in D$. Let $a \in D$ be limit of a sequence $\{a_n\}$. Then*

$$\lim_{n\to\infty}(a_n d) = ad$$

$$\lim_{n\to\infty}(d a_n) = da$$

Proof. Statement of the theorem is trivial, however I give this proof for completeness sake. Since $a \in D$ is limit of the sequence $\{a_n\}$, then according to definition 2.1.3 for given $\epsilon \in R$, $\epsilon > 0$, there exists positive integer $n_0$ such, that

$$|a_n - a| < \frac{\epsilon}{|d|}$$

for every $n > n_0$. According to definition 2.1.2 the statement of theorem follows from inequalities

$$|a_n d - ad| = |(a_n - a)d| = |a_n - a||d| < \frac{\epsilon}{|d|}|d| = \epsilon$$

$$|d a_n - da| = |d(a_n - a)| = |d||a_n - a| < |d|\frac{\epsilon}{|d|} = \epsilon$$

for any $n > n_0$. □

Definition 2.1.5. *Let $D$ be normed ring. The sequence $\{a_n\}$, $a_n \in D$ is called* **fundamental** *or* **Cauchy sequence***, if for every $\epsilon \in R$, $\epsilon > 0$ there exists positive integer $n_0$ depending on $\epsilon$ and such, that $|a_p - a_q| < \epsilon$ for every $p$, $q > n_0$.* □

Definition 2.1.6. *Normed ring $D$ is called* **complete** *if any fundamental sequence of elements of ring $D$ converges, i.e. has limit in ring $D$.* □

Later on, speaking about normed ring of characteristic 0, we will assume that homeomorphism of field of rational numbers $Q$ into ring $D$ is defined.

Theorem 2.1.7. *Complete ring $D$ of characteristic 0 contains as subfield an isomorphic image of the field $R$ of real numbers. It is customary to identify it with $R$.*

Proof. Consider fundamental sequence of rational numbers $\{p_n\}$. Let $p'$ be limit of this sequence in ring $D$. Let $p$ be limit of this sequence in field $R$. Since immersion of field $Q$ into division ring $D$ is homeomorphism, then we may identify $p' \in D$ and $p \in R$. □

Theorem 2.1.8. *Let $D$ be complete ring of characteristic 0 and let $d \in D$. Then any real number $p \in R$ commute with $d$.*

Proof. Let us represent real number $p \in R$ as fundamental sequence of rational numbers $\{p_n\}$. Statement of theorem follows from chain of equations

$$pd = \lim_{n\to\infty}(p_n d) = \lim_{n\to\infty}(d p_n) = dp$$

based on statement of theorem 2.1.4. □

## 2.2. Topological $D$-Algebra

Definition 2.2.1. *Given a topological commutative ring $D$ and $D$-algebra $A$ such that $A$ has a topology compatible with the structure of the additive group of $A$ and maps*

$$(a, v) \in D \times A \to av \in A$$
$$(v, w) \in A \times A \to vw \in A$$

*are continuous, then $\overline{V}$ is called a* **topological $D$-algebra**[2.3]. □

Definition 2.2.2. **Norm on $D$-algebra** *$A$ over normed commutative ring $D$*[2.4] *is a map*

$$a \in A \to |a| \in R$$

*which satisfies the following axioms*

- $|a| \geq 0$
- $|a| = 0$ *if, and only if,* $a = 0$
- $|a + b| \leq |a| + |b|$
- $|ab| = |a|\ |b|$
- $|da| = |d|\ |a|$, $d \in D$, $a \in A$

*If $D$ is a normed commutative ring, $D$-algebra $A$, endowed with the structure defined by a given norm on $A$, is called* **normed $D$-algebra**. □

Definition 2.2.3. *Let $A$ be normed $D$-algebra. Element $a \in A$ is called* **limit of a sequence** $\{a_n\}$

$$a = \lim_{n\to\infty} a_n$$

*if for every $\epsilon \in R$, $\epsilon > 0$ there exists positive integer $n_0$ depending on $\epsilon$ and such, that $|a_n - a| < \epsilon$ for every $n > n_0$.* □

Definition 2.2.4. *Let $A$ be normed $D$-algebra. The sequence $\{a_n\}$, $a_n \in A$, is called* **fundamental** *or* **Cauchy sequence**, *if for every $\epsilon \in R$, $\epsilon > 0$, there exists positive integer $n_0$ depending on $\epsilon$ and such, that $|a_p - a_q| < \epsilon$ for every $p$, $q > n_0$.* □

Definition 2.2.5. *Normed $D$-algebra $A$ is called* **Banach $D$-algebra** *if any fundamental sequence of elements of algebra $A$ converges, i.e. has limit in algebra $A$.* □

Definition 2.2.6. *Let $A$ be Banach $D$-algebra. Set of elements $a \in A$, $|a| = 1$, is called* **unit sphere in algebra $A$**. □

Definition 2.2.7. *A map*

$$f : A_1 \to A_2$$

*of Banach $D_1$-algebra $A_1$ with norm $|x|_1$ into Banach $D_2$-algebra $A_2$ with norm $|y|_2$ is called* **continuous**, *if for every as small as we please $\epsilon > 0$ there exist such $\delta > 0$, that*

$$|x' - x|_1 < \delta$$

*implies*

$$|f(x') - f(x)|_2 < \epsilon$$

[2.3] I made definition according to definition from [3], p. TVS I.1

[2.4] I made definition according to definition from [2], IX, §3.3

□

Definition 2.2.8. *Let*

$$f : A_1 \to A_2$$

*be map of Banach $D_1$-algebra $A_1$ with norm $|x|_1$ into Banach $D_2$-algebra $A_2$ with norm $|y|_2$. Value*

$$\|f\| = sup\frac{|f(x)|_2}{|x|_1} \tag{2.2.1}$$

*is called* **norm of map** *$f$.* □

Theorem 2.2.9. *Let*

$$f : A_1 \to A_2$$

*be linear map of Banach $D_1$-algebra $A_1$ with norm $|x|_1$ into Banach $D_2$-algebra $A_2$ with norm $|y|_2$. Then*

$$\|f\| = sup\{|f(x)|_2 : |x|_1 = 1\} \tag{2.2.2}$$

Proof. From definitions [1]-4.2.1, [1]-6.1.1 and theorems 2.1.7, 2.1.8, it follows that

$$f(rx) = rf(x) \quad r \in R \tag{2.2.3}$$

From the equation (2.2.3) and the definition 2.2.2 it follows that

$$\frac{|f(rx)|_2}{|rx|_1} = \frac{|r|\ |f(x)|_2}{|r|\ |x|_1} = \frac{|f(x)|_2}{|x|_1}$$

Assuming $r = \frac{1}{|x|_1}$, we get

$$\frac{|f(x)|_2}{|x|_1} = \left|f\left(\frac{x}{|x|_1}\right)\right|_2 \tag{2.2.4}$$

Equation (2.2.2) follows from equations (2.2.4) and (2.2.1). □

Theorem 2.2.10. *Let*

$$f : A_1 \to A_2$$

*be linear map of Banach $D_1$-algebra $A_1$ with norm $|x|_1$ into Banach $D_2$-algebra $A_2$ with norm $|y|_2$. Since $\|f\| < \infty$, then map $f$ is continuous.*

Proof. Since map $f$ is linear, then according to definition 2.2.8

$$|f(x) - f(y)|_2 = |f(x - y)|_2 \le \|f\|\ |x - y|_1$$

Let us assume arbitrary $\epsilon > 0$. Assume $\delta = \frac{\epsilon}{\|f\|}$. Then

$$|f(x) - f(y)|_2 \le \|f\|\ \delta = \epsilon$$

follows from inequality

$$|x - y|_1 < \delta$$

According to definition 2.2.7 map $f$ is continuous. □

## 2.3. The Derivative of Map in Algebra

Definition 2.3.1. *Let $A$ be Banach $D$-algebra. The map*

$$f : A \to A$$

*is called* **differentiable in the Gâteaux sense** *on the set $U \subset A$, if at every point $x \in U$ the increment of the map $f$ can be represented as*

$$f(x+a) - f(x) = \partial f(x) \circ a + o(a) = \frac{\partial f(x)}{\partial x} \circ a + o(a) \tag{2.3.1}$$

*where* **the Gâteaux derivative $\partial f(x)$ of map** $f$ *is linear map of increment $a$ and $o : A \to A$ is such continuous map that*

$$\lim_{a \to 0} \frac{|o(a)|}{|a|} = 0$$

□

Remark 2.3.2. *According to definition 2.3.1 for given $x$, the Gâteaux derivative $\partial f(x) \in \mathcal{L}(A; A)$. Therefore, the Gâteaux derivative of map $f$ is map*

$$\partial f : A \to \mathcal{L}(A; A)$$

*Expressions $\partial f(x)$ and $\dfrac{\partial f(x)}{\partial x}$ are different notations for the same map. We will use notation $\dfrac{\partial f(x)}{\partial x}$ to underline that this is the Gâteaux derivative with respect to variable $x$.* □

Theorem 2.3.3. *It is possible to represent the Gâteaux derivative of map $f$ as*[2.5]

$$\begin{aligned} \frac{\partial f(x)}{\partial x} \circ a &= \left( \frac{\partial_{s \cdot 0} f(x)}{\partial x} \otimes \frac{\partial_{s \cdot 1} f(x)}{\partial x} \right) \circ a \\ &= \frac{\partial_{s \cdot 0} f(x)}{\partial x} a \frac{\partial_{s \cdot 1} f(x)}{\partial x} \end{aligned} \tag{2.3.3}$$

*Expression $\dfrac{\partial_{s \cdot p} f(x)}{\partial x}$, $p = 0, 1$, is called* **component of the Gâteaux derivative** *of map $f(x)$.*

Proof. The theorem follows from the definitions 2.3.1 and from the theorem [1]-6.4.5. □

From definitions [1]-4.2.1, [1]-6.1.1 , 2.3.1 and the theorem 2.1.7 it follows

$$\partial f(x) \circ (ra) = r \partial f(x) \circ a \tag{2.3.4}$$

$$r \in R \quad r \neq 0 \quad a \in A \quad a \neq 0$$

---

[2.5] Formally, we have to write the differential of the map in the form

$$\frac{\partial f(x)}{\partial x} \circ a = \left( \frac{\partial_{k \cdot s \cdot 0} f(x)}{\partial x} \otimes \frac{\partial_{k \cdot s \cdot 1} f(x)}{\partial x} \right) \circ I_k \circ a = \frac{\partial_{k \cdot s \cdot 0} f(x)}{\partial x} (I_k \circ a) \frac{\partial_{k \cdot s \cdot 1} f(x)}{\partial x} \tag{2.3.2}$$

However, for instance, in the theory of functions of complex variable we consider only linear maps generated by map $I_0 \circ z = z$. Therefore, exploring derivatives, we also restrict ourselves to linear maps generated by the map $I_0$. To write expressions in the general case is not difficult.

Combining equation (2.3.4) and definition 2.3.1, we get known definition of the Gâteaux derivative

$$\partial f(x)\circ a=\lim_{t\to 0,\ t\in R}(t^{-1}(f(x+ta)-f(x))) \tag{2.3.5}$$

Definitions of the Gâteaux derivative (2.3.1) and (2.3.5) are equivalent. Using this equivalence we tell that map $f$ is called differentiable in the Gâteaux sense on the set $U\subset D$, if at every point $x\in U$ the increment of the map $f$ can be represented as

$$f(x+ta)-f(x)=t\partial f(x)\circ a+o(t) \tag{2.3.6}$$

where $o:R\to A$ is such continuous map that

$$\lim_{t\to 0}\frac{|o(t)|}{|t|}=0$$

Since infinitesimal $a$ in the equation (2.3.1) is differential $dx$, then equation (2.3.1) becomes definition of **the Gâteaux differential**

$$\begin{aligned}\frac{\partial f(x)}{\partial x}\circ dx&=\left(\frac{\partial_{s\cdot 0}f(x)}{\partial x}\otimes\frac{\partial_{s\cdot 1}f(x)}{\partial x}\right)\circ dx\\&=\frac{\partial_{s\cdot 0}f(x)}{\partial x}dx\frac{\partial_{s\cdot 1}f(x)}{\partial x}\end{aligned} \tag{2.3.7}$$

THEOREM 2.3.4. *Let $A$ be Banach $D$-algebra. Let $\overline{\overline{e}}$ be basis of algebra $A$ over ring $D$.* **Standard representation of the Gâteaux derivative of mapping**

$$f:A\to A$$

*has form*

$$\frac{\partial f(x)}{\partial x}=\frac{\partial^{ij}f(x)}{\partial x}\overline{e}_i\otimes\overline{e}_j \tag{2.3.8}$$

*Expression* $\frac{\partial^{ij}f(x)}{\partial x}$ *in equation (2.3.8) is called* **standard component of the Gâteaux derivative** *of mapping $f$.*

PROOF. Statement of theorem folows from the statement [1]-6.4.5.2. □

THEOREM 2.3.5. *Let $A$ be Banach $D$-algebra. Let $\overline{\overline{e}}$ be basis of algebra $A$ over ring $D$. Then it is possible to represent the Gâteaux derivative of map*

$$f:D\to D$$

*as*

$$\frac{\partial f(x)}{\partial x}\circ dx=dx^i\frac{\partial f^j}{\partial x^i}\overline{e}_j \tag{2.3.9}$$

*where $dx\in A$ has expansion*

$$dx=dx^i\ e_i\qquad dx^i\in D$$

*relative to basis $\overline{\overline{e}}$ and Jacobian matrix of map $f$ has form*

$$\frac{\partial f^j}{\partial x^i}=\frac{\partial^{kr}f(x)}{\partial x}C^p_{ki}C^j_{pr} \tag{2.3.10}$$

PROOF. Statement of theorem follows from the theorem [1]-6.4.5. □

Theorem 2.3.6. *Let $A$ be Banach $D$-algebra. Let $f$, $g$ be differentiable maps*

$$f : A \to A \quad g : A \to A$$

*The map*

$$f + g : A \to A$$

*is differentiable and the Gâteaux derivative satisfies to relationship*

$$\partial(f+g)(x) = \partial f(x) + \partial g(x) \tag{2.3.11}$$

Proof. According to the definition (2.3.5),

$$\begin{aligned}
\partial(f+g)(x)\circ a &= \lim_{t\to 0,\ t\in R}(t^{-1}((f+g)(x+ta)-(f+g)(x)))\\
&= \lim_{t\to 0,\ t\in R}(t^{-1}(f(x+ta)+g(x+ta)-f(x)-g(x)))\\
&= \lim_{t\to 0,\ t\in R}(t^{-1}(f(x+ta)-f(x)))\\
&+ \lim_{t\to 0,\ t\in R}(t^{-1}(g(x+ta)-g(x)))\\
&= \partial f(x)\circ a + \partial g(x)\circ a
\end{aligned} \tag{2.3.12}$$

The equation (2.3.11) follows from the equation (2.3.12). □

Theorem 2.3.7. *Let $A$ be Banach $D$-algebra. Let*

$$h : A \times A \to A$$

*be continuous bilinear map. Let $f$, $g$ be differentiable maps*

$$f : A \to A \quad g : A \to A$$

*The map*

$$h(f, g) : A \to A$$

*is differentiable and the Gâteaux differential satisfies to relationship*

$$\partial h(f(x), g(x)) \circ dx = h(\partial f(x) \circ dx, g(x)) + h(f(x), \partial g(x) \circ dx) \tag{2.3.13}$$

Proof. Equation (2.3.13) follows from chain of equations

$$\begin{aligned}
\partial h(f(x), g(x)) \circ a &= \lim_{t\to 0}(t^{-1}(h(f(x+ta), g(x+ta)) - h(f(x), g(x))))\\
&= \lim_{t\to 0}(t^{-1}(h(f(x+ta), g(x+ta)) - h(f(x), g(x+ta))))\\
&+ \lim_{t\to 0}(t^{-1}(h(f(x), g(x+ta)) - h(f(x), g(x))))\\
&= h(\lim_{t\to 0} t^{-1}(f(x+ta) - f(x)), g(x))\\
&+ h(f(x), \lim_{t\to 0} t^{-1}(g(x+ta) - g(x)))
\end{aligned}$$

based on definition (2.3.5). □

Convention 2.3.8. *Given bilinear map*

$$h : A \times A \to A$$

*we consider following maps*

$$h_1 : \mathcal{L}(D; A; A) \times A \to \mathcal{L}(D; A; A)$$
$$h_2 : A \times \mathcal{L}(D; A; A) \to \mathcal{L}(D; A; A)$$

*defined by equation*

$$h_1(f,v)\circ u = h(f\circ u, v)$$
$$h_2(u,f)\circ v = h(u, f\circ v)$$

*We will use letter $h$ to denote maps $h_1$, $h_2$.* □

THEOREM 2.3.9. *Let $A$ be Banach $D$-algebra. Let*

$$h : A\times A\rightarrow A$$

*be continuous bilinear map. Let $f$, $g$ be differentiable maps*

$$f:A\rightarrow A \quad g:A\rightarrow A$$

*The map*

$$h(f,g):A\rightarrow A$$

*is differentiable and the Gâteaux derivative satisfies to relationship*

$$\partial h(f(x),g(x)) = h(\partial f(x),g(x)) + h(f(x),\partial g(x)) \tag{2.3.14}$$

PROOF. Equation (2.3.14) follows from the equation (2.3.13) and from the convention 2.3.8. □

THEOREM 2.3.10. *Let $A$ be Banach $D$-algebra. Let $f$, $g$ be differentiable maps*

$$f:A\rightarrow A \quad g:A\rightarrow A$$

*The Gâteaux derivative satisfies to relationship*

$$\frac{\partial f(x)g(x)}{\partial x}\circ dx = \left(\frac{\partial f(x)}{\partial x}\circ dx\right)\ g(x) + f(x)\ \left(\frac{\partial g(x)}{\partial x}\circ dx\right) \tag{2.3.15}$$

$$\frac{\partial f(x)g(x)}{\partial x} = \frac{\partial f(x)}{\partial x}\ g(x) + f(x)\ \frac{\partial g(x)}{\partial x} \tag{2.3.16}$$

PROOF. The theorem follows from theorems 2.3.7, 2.3.9 and the definition [1]-5.1.1. □

THEOREM 2.3.11. *Let $A$ be Banach $D$-algebra. Let the Gâteaux derivative of map*

$$f:A\rightarrow A$$

*have expansion*

$$\frac{\partial f(x)}{\partial x} = \frac{\partial_{s\cdot 0} f(x)}{\partial x}\otimes\frac{\partial_{s\cdot 1} f(x)}{\partial x} \tag{2.3.17}$$

*Let the Gâteaux derivative of map*

$$g:A\rightarrow A$$

*have expansion*

$$\frac{\partial g(x)}{\partial x} = \frac{\partial_{t\cdot 0} g(x)}{\partial x}\otimes\frac{\partial_{t\cdot 1} g(x)}{\partial x} \tag{2.3.18}$$

*The Gâteaux derivative of map $f(x)g(x)$ have form*

$$\frac{\partial f(x)g(x)}{\partial x} = \frac{\partial_{s\cdot 0} f(x)}{\partial x}\otimes\left(\frac{\partial_{s\cdot 1} f(x)}{\partial x}g(x)\right) + \left(f(x)\frac{\partial_{t\cdot 0} g(x)}{\partial x}\right)\otimes\frac{\partial_{t\cdot 1} g(x)}{\partial x} \tag{2.3.19}$$

$$\begin{aligned}&\frac{\partial_{s\cdot 0}f(x)g(x)}{\partial x}=\frac{\partial_{s\cdot 0}f(x)}{\partial x} && \frac{\partial_{t\cdot 0}f(x)g(x)}{\partial x}=f(x)\frac{\partial_{t\cdot 0}g(x)}{\partial x}\\ &\frac{\partial_{s\cdot 1}f(x)g(x)}{\partial x}=\frac{\partial_{s\cdot 1}f(x)}{\partial x}g(x) && \frac{\partial_{t\cdot 1}f(x)g(x)}{\partial x}=\frac{\partial_{t\cdot 1}g(x)}{\partial x}\end{aligned}\tag{2.3.20}$$

Proof. Let us substitute (2.3.17) and (2.3.18) into equation (2.3.16)

$$\begin{aligned}\frac{\partial f(x)g(x)}{\partial x}\circ a&=\left(\frac{\partial f(x)}{\partial x}\circ a\right)\ g(x)+f(x)\ \left(\frac{\partial g(x)}{\partial x}\circ a\right)\\ &=\frac{\partial_{s\cdot 0}f(x)}{\partial x}a\frac{\partial_{s\cdot 1}f(x)}{\partial x}g(x)+f(x)\frac{\partial_{t\cdot 0}g(x)}{\partial x}a\frac{\partial_{t\cdot 1}g(x)}{\partial x}\end{aligned}\tag{2.3.21}$$

Based (2.3.21), we define equations (2.3.20). □

Theorem 2.3.12. *Let $A$ be Banach $D$-algebra. Let $f$, $g$ be differentiable maps*

$$f:A\to A\quad g:A\to A$$

*The Gâteaux derivative satisfies to relationship*

$$\frac{\partial f(x)\otimes g(x)}{\partial x}\circ dx=\left(\frac{\partial f(x)}{\partial x}\circ dx\right)\otimes g(x)+f(x)\otimes\left(\frac{\partial g(x)}{\partial x}\circ dx\right)\tag{2.3.22}$$

$$\frac{\partial f(x)\otimes g(x)}{\partial x}=\frac{\partial f(x)}{\partial x}\otimes g(x)+f(x)\otimes\frac{\partial g(x)}{\partial x}\tag{2.3.23}$$

Proof. The theorem follows from the theorems 2.3.7, 2.3.9, [1]-4.4.5 and the definition [1]-5.1.1. □

Remark 2.3.13. *Let*

$$\frac{\partial f(x)}{\partial x}=\frac{\partial_{s\cdot 0}f(x)}{\partial x}\otimes\frac{\partial_{s\cdot 1}f(x)}{\partial x}\tag{2.3.24}$$

$$\frac{\partial g(x)}{\partial x}=\frac{\partial_{t\cdot 0}g(x)}{\partial x}\otimes\frac{\partial_{t\cdot 1}g(x)}{\partial x}\tag{2.3.25}$$

*Then*

$$\frac{\partial f(x)\otimes g(x)}{\partial x}=\frac{\partial_{s\cdot 0}f(x)}{\partial x}\otimes\frac{\partial_{s\cdot 1}f(x)}{\partial x}\otimes g(x)+f(x)\otimes\frac{\partial_{t\cdot 0}g(x)}{\partial x}\otimes\frac{\partial_{t\cdot 1}g(x)}{\partial x}\tag{2.3.26}$$

*We do not write brackets, because tensor product is associative and distributive over addition (theorems [1]-3.3.11, [1]-3.4.5 ).* □

Theorem 2.3.14. *Let $A$ be Banach $D$-algebra. If the Gâteaux derivative $\partial f(x)$ exists in point $x$ and has finite norm, then map $f$ is continuous at point $x$.*

Proof. From definition 2.2.8 it follows

$$|\partial f(x)\circ a|\le\|\partial f(x)\|\ |a|\tag{2.3.27}$$

From (2.3.1), (2.3.27) it follows

$$|f(x+a)-f(x)|<|a|\ \|\partial f(x)\|\tag{2.3.28}$$

Let us assume arbitrary $\epsilon>0$. Assume

$$\delta=\frac{\epsilon}{\|\partial f(x)\|}$$

Then from inequality

$$|a|<\delta$$

it follows

$$|f(x+a)-f(x)| \le \|\partial f(x)\|\ \delta = \epsilon$$

According to definition 2.2.7 map $f$ is continuous at point $x$. □

THEOREM 2.3.15. *Let $A$ be Banach $D$-algebra. Let map*

$$f: A \to A$$

*be differentiable in the Gâteaux sense at point $x$. Then*

$$\partial f(x) \circ 0 = 0$$

PROOF. The theorem follows from the definitions 2.3.1 and from the theorem [1]-4.2.5. □

THEOREM 2.3.16. *Let $A$ be Banach $D$-algebra. Let map*

$$f: A \to A$$

*be differentiable in the Gâteaux sense at point $x$ and norm of the Gâteaux derivative of map $f$ be finite*

(2.3.29) $$\|\partial f(x)\| = F \le \infty$$

*Let map*

$$g: A \to A$$

*be differentiable in the Gâteaux sense at point*

(2.3.30) $$y = f(x)$$

*and norm of the Gâteaux derivative of map $g$ be finite*

(2.3.31) $$\|\partial g(y)\| = G \le \infty$$

*Map*

$$(g \circ f)(x) = g(f(x))$$

*is differentiable in the Gâteaux sense at point $x$*

(2.3.32) $$\begin{cases} \partial(g\circ f)(x) = \partial g(y) \circ \partial f(x) \\ \partial(g\circ f)(x)\circ a = \partial g(y)\circ \partial f(x) \circ a \end{cases}$$

(2.3.33) $$\begin{cases} \dfrac{\partial_{st\cdot 0}(g\circ f)(x)}{\partial x} = \dfrac{\partial_{s\cdot 0} g(f(x))}{\partial f(x)} \dfrac{\partial_{t\cdot 0} f(x)}{\partial x} \\ \dfrac{\partial_{st\cdot 1}(g\circ f)(x)}{\partial x} = \dfrac{\partial_{t\cdot 1} f(x)}{\partial x} \dfrac{\partial_{s\cdot 1} g(f(x))}{\partial f(x)} \end{cases}$$

PROOF. According to definition 2.3.1

(2.3.34) $$g(y+b) - g(y) = \partial g(y)\circ b + o_1(b)$$

where $o_1: A \to A$ is such continuous map that

$$\lim_{b\to 0} \frac{|o_1(b)|}{|b|} = 0$$

According to definition 2.3.1

(2.3.35) $$f(x+a) - f(x) = \partial f(x)\circ a + o_2(a)$$

where $o_2 : A \to A$ is such continuous map that

$$\lim_{a\to 0}\frac{|o_2(a)|}{|a|}=0$$

According to (2.3.35) increment $a$ of value $x \in A$ leads to increment

$$b=\partial f(x)\circ a+o_2(a) \tag{2.3.36}$$

of value $y$. Using (2.3.30), (2.3.36) in equation (2.3.34), we get

$$\begin{aligned}&g(f(x+a))-g(f(x))\\=&g(f(x)+\partial f(x)\circ a+o_2(a))-g(f(x))\\=&\partial g(f(x))\circ(\partial f(x)\circ a+o_2(a))-o_1(\partial f(x)\circ a+o_2(a))\end{aligned} \tag{2.3.37}$$

According to definitions 2.3.1, [1]-4.2.1, [1]-6.1.1 , from equation (2.3.37) it follows

$$\begin{aligned}&g(f(x+a))-g(f(x))\\=&\partial g(f(x))\circ\partial f(x)\circ a+\partial g(f(x))\circ o_2(a)-o_1(\partial f(x)\circ a+o_2(a))\end{aligned} \tag{2.3.38}$$

According to definition 2.2.2

$$\begin{aligned}&\lim_{a\to 0}\frac{|\partial g(f(x))\circ o_2(a)-o_1(\partial f(x)\circ a+o_2(a))|}{|a|}\\\le&\lim_{a\to 0}\frac{|\partial g(f(x))\circ o_2(a)|}{|a|}+\lim_{a\to 0}\frac{|o_1(\partial f(x)\circ a+o_2(a))|}{|a|}\end{aligned} \tag{2.3.39}$$

From (2.3.31) it follows that

$$\lim_{a\to 0}\frac{|\partial g(f(x))\circ o_2(a)|}{|a|}\le G\lim_{a\to 0}\frac{|o_2(a)|}{|a|}=0 \tag{2.3.40}$$

From (2.3.29) it follows that

$$\begin{aligned}&\lim_{a\to 0}\frac{|o_1(\partial f(x)\circ a+o_2(a))|}{|a|}\\=&\lim_{a\to 0}\frac{|o_1(\partial f(x)\circ a+o_2(a))|}{|\partial f(x)\circ a+o_2(a)|}\lim_{a\to 0}\frac{|\partial f(x)\circ a+o_2(a)|}{|a|}\\\le&\lim_{a\to 0}\frac{|o_1(\partial f(x)\circ a+o_2(a))|}{|\partial f(x)\circ a+o_2(a)|_2}\lim_{a\to 0}\frac{\|\partial f(x)\||a|+|o_2(a)|}{|a|}\\=&\lim_{a\to 0}\frac{|o_1(\partial f(x)\circ a+o_2(a))|}{|\partial f(x)\circ a+o_2(a)|}\|\partial f(x)\|\end{aligned}$$

According to the theorem 2.3.15

$$\lim_{a\to 0}(\partial f(x)\circ a)+o_2(a))=0$$

Therefore,

$$\lim_{a\to 0}\frac{|o_1(\partial f(x)\circ a+o_2(a))|}{|a|}=0 \tag{2.3.41}$$

From equations (2.3.39), (2.3.40), (2.3.41) it follows

$$\lim_{a\to 0}\frac{|\partial g(f(x))\circ o_2(a)-o_1(\partial f(x)\circ a+o_2(a))|}{|a|}=0 \tag{2.3.42}$$

According to definition 2.3.1

$$(g \circ f)(x+a) - (g \circ f)(x) = \partial(g \circ f)(x) \circ a + o(a) \tag{2.3.43}$$

where $o : A \to A$ is such continuous map that

$$\lim_{a \to 0} \frac{|o(a)|}{|a|} = 0$$

Equation (2.3.32) follows from (2.3.38), (2.3.42), (2.3.43).

From equation (2.3.32) and theorem 2.3.3 it follows that

$$\begin{aligned}
&\frac{\partial_{st\cdot 0}(g \circ f)(x)}{\partial x} \ a \ \frac{\partial_{st\cdot 1}(g \circ f)(x)}{\partial x} \\
&= \frac{\partial_{s\cdot 0} g(f(x))}{\partial f(x)} \ (\partial f(x) \circ a) \ \frac{\partial_{s\cdot 1} g(f(x))}{\partial f(x)} \\
&= \frac{\partial_{s\cdot 0} g(f(x))}{\partial f(x)} \frac{\partial_{t\cdot 0} f(x)}{\partial x} \ a \ \frac{\partial_{t\cdot 1} f(x)}{\partial x} \ \frac{\partial_{s\cdot 1} g(f(x))}{\partial f(x)}
\end{aligned} \tag{2.3.44}$$

(2.3.33) follow from equation (2.3.44). □

CHAPTER 3

# Derivative of Second Order of Map of $D$-Algebra

## 3.1. Derivative of Second Order of Map of $D$-Algebra

Let $D$ be the complete commutative ring of characteristic 0. Let $A$ be associative $D$-algebra. Let

$$f:A\rightarrow A$$

function differentiable in the Gâteaux sense. According to remark 2.3.2 the Gâteaux derivative is map

$$\partial f:A\rightarrow\mathcal{L}(A;A)$$

According to the theorem [1]-6.2.5 and the definition 2.2.8, set $\mathcal{L}(A;A)$ is Banach $D$-algebra. Therefore, we may consider the question, if map $\partial f$ is differentiable in the Gâteaux sense.

According to definition 2.3.1

$$(\partial f\circ(x+a_2))\circ a_1-(\partial f\circ x)\circ a_1=\partial(\partial f(x)\circ a_1)\circ a_2+o_2(a_2) \tag{3.1.1}$$

where $o_2:A\rightarrow\mathcal{L}(A;A)$ is such continuous map, that

$$\lim_{a_2\rightarrow 0}\frac{\|o_2(a_2)\|}{|a_2|}=0$$

According to definition 2.3.1 the mapping $\partial(\partial f(x)\circ a_1)\circ a_2$ is linear map of variable $a_2$. From equation (3.1.1) it follows that mapping $\partial(\partial f(x)\circ a_1)\circ a_2$ is linear mapping of variable $a_1$. Therefore, the mapping $\partial(\partial f(x)\circ a_1)\circ a_2$ is bilinear mapping.

Definition 3.1.1. *Polylinear mapping*

$$\partial^2 f(x)\circ(a_1;a_2)=\frac{\partial^2 f(x)}{\partial x^2}\circ(a_1;a_2)=\partial(\partial f(x)\circ a_1)\circ a_2 \tag{3.1.2}$$

*is called* **the Gâteaux derivative of second order** *of map $f$.* □

Remark 3.1.2. *According to definition 3.1.1 for given $x$ the Gâteaux derivative of second order $\partial^2 f(x)\in\mathcal{L}(A,A;A)$. Therefore, the Gâteaux derivative of second order of map $f$ is mapping*

$$\partial^2 f:A\rightarrow\mathcal{L}(A,A;A)$$

*According to the theorem [1]-4.4.4, we may consider also expression*

$$\partial^2 f(x)\circ(a_1\otimes a_2)=\partial^2 f(x)\circ(a_1;a_2)$$

*Then*

$$\partial^2 f(x)\in\mathcal{L}(A\otimes A;A)$$
$$\partial^2 f:A\rightarrow\mathcal{L}(A\otimes A;A)$$

*We use the same notation for mapping because of the nature of the argument it is clear what kind of mapping we consider.* □

Theorem 3.1.3. *It is possible to represent the Gâteaux derivative of second order of map $f$ as*

$$\partial^2 f(x)\circ(a_1;a_2)=\left(\frac{\partial^2_{s\cdot 0}f(x)}{\partial x^2}\otimes\frac{\partial^2_{s\cdot 1}f(x)}{\partial x^2}\otimes\frac{\partial^2_{s\cdot 2}f(x)}{\partial x^2},\sigma_s\right)\circ(a_1;a_2)$$
$$=\frac{\partial^2_{s\cdot 0}f(x)}{\partial x^2}\sigma_s(a_1)\frac{\partial^2_{s\cdot 1}f(x)}{\partial x^2}\sigma_s(a_2)\frac{\partial^2_{s\cdot 2}f(x)}{\partial x^2}$$

*Expression* [3.1]

$$\frac{\partial^2_{s\cdot p}f(x)}{\partial x^2}\quad p=0,1,2$$

*is called* **component of the Gâteaux derivative of second order** *of map $f(x)$.*

Proof. Corollary of definition 3.1.1 and theorem [1]-6.6.6. □

By induction, assuming that we defined the Gâteaux derivative $\partial^{n-1}f(x)$ of order $n-1$, we define

$$\partial^n f(x)\circ(a_1;...;a_n)=\frac{\partial^n f(x)}{\partial x^n}\circ(a_1;...;a_n)\tag{3.1.3}$$
$$=\ \partial(\partial^{n-1}f(x)\circ(a_1;...;a_{n-1}))\circ a_n$$

**the Gâteaux derivative of order** $n$ of map $f$. We also assume $\partial^0 f(x)=f(x)$.

## 3.2. Taylor Series

Let $D$ be the complete commutative ring of characteristic 0. Let $A$ be associative $D$-algebra. Let $p_k(x)$ be the monomial of power $k$, $k>0$, in one variable over $D$-algebra $A$.

It is evident that monomial of power 0 has form $a_0$, $a_0\in A$. For $k>0$,

$$p_k(x)=p_{k-1}(x)xa_k$$

where $a_k\in A$. Actually, last factor of monomial $p_k(x)$ is either $a_k\in A$, or has form $x^l$, $l\ge 1$. In the later case we assume $a_k=1$. Factor preceding $a_k$ has form $x^l$, $l\ge 1$. We can represent this factor as $x^{l-1}x$. Therefore, we proved the statement.

In particular, monomial of power 1 has form $p_1(x)=a_0xa_1$.

Without loss of generality, we assume $k=n$.

Theorem 3.2.1. *For any $m>0$ the following equation is true*

$$\begin{aligned}\partial^m(f(x)x)\circ(h_1;...;h_m)=&\ \partial^m f(x)\circ(h_1;...;h_m)x\\&+\partial^{m-1}f(x)\circ(\widehat{h_1};...;h_{m-1};h_m)h_1+...\\&+\partial^{m-1}f(x)\circ(h_1;...;h_{m-1};\widehat{h_m})h_m\end{aligned}\tag{3.2.1}$$

*where symbol $\widehat{h^i}$ means absense of variable $h^i$ in the list.*

[3.1] We suppose

$$\frac{\partial^2_{s\cdot p}f(x)}{\partial x^2}=\frac{\partial^2_{s\cdot p}f(x)}{\partial x\partial x}$$

PROOF. For $m=1$, this is corollary of equation (2.3.16)

$$\partial(f(x)x)\circ h_1=(\partial f(x)\circ h_1)x+f(x)h_1$$

Assume, (3.2.1) is true for $m-1$. Then

$$\begin{aligned}\partial^{m-1}(f(x)x)\circ(h_1;...;h_{m-1})=&\partial^{m-1}f(x)\circ(h_1;...;h_{m-1})x\\&+\partial^{m-2}f(x)\circ(\widehat{h_1};...;h_{m-2};h_{m-1})h_1+...\\&+\partial^{m-2}f(x)\circ(h_1;...;h_{m-2};\widehat{h_{m-1}})h_{m-1}\end{aligned}$$

Since $\partial h_i=0$, then, using the equation (2.3.16), we get
(3.2.2)

$$\begin{aligned}\partial^{m}(f(x)x)\circ(h_1;...;h_{m-1};h_m)=&\partial^{m}f(x)\circ(h_1;...;h_{m-1};h_m)x\\&+\partial^{m-1}f(x)\circ(h_1;...;h_{m-2};h_{m-1})h_m\\&+\partial^{m-1}f(x)\circ(\widehat{h_1};...;h_{m-2};h_{m-1};h_m)h_1+...\\&+\partial^{m-1}f(x)\circ(h_1;...;h_{m-2};\widehat{h_{m-1}};h_m)h_{m-1}\end{aligned}$$

The difference between equations (3.2.1) and (3.2.2) is only in form of presentation. We proved the theorem. □

THEOREM 3.2.2. *For any $n\ge 0$ following equation is true*

$$\partial^{n+1}p_n(x)=0$$

PROOF. Since $p_0(x)=a_0$, $a_0\in D$, then for $n=0$ theorem is corollary of definition 2.3.1. Let statement of theorem is true for $n-1$. According to theorem 3.2.1, when $f(x)=p_{n-1}(x)$, we get

$$\begin{aligned}\partial^{n+1}p_n(x)(h_1;...;h_{n+1})=&\partial^{n+1}(p_{n-1}(x)xa_n)(h_1;...;h_{n+1})\\=&\partial^{n+1}p_{n-1}(x)(h_1;...;h_{n+1})xa_n\\&+\partial^{n}p_{n-1}(x)(\widehat{h_1};...;h_n;h_{n+1})h_1a_n+...\\&+\partial^{n}p_{n-1}(x)(h_1;...;h_n;\widehat{h_{n+1}})h_{n+1}a_n\end{aligned}$$

According to suggestion of induction all monomials are equal 0. □

THEOREM 3.2.3. *If $m<n$, then following equation is true*

$$\partial^m p_n(0)=0$$

PROOF. For $n=1$ following equation is true

$$\partial^0 p_1(0)=a_0xa_1=0$$

Assume that statement is true for $n-1$. Then according to theorem 3.2.1

$$\begin{aligned}\partial^m(p_{n-1}(x)xa_n)(h_1;...;h_m)=&\partial^m p_{n-1}(x)(h_1;...;h_m)xa_n\\&+\partial^{m-1}p_{n-1}(x)(\widehat{h_1};...;h_{m-1};h_m)h_1a_n+...\\&+\partial^{m-1}p_{n-1}(x)(h_1;...;h_{m-1};\widehat{h_m})h_ma_n\end{aligned}$$

First term equal 0 because $x=0$. Because $m-1<n-1$, then rest terms equal 0 according to suggestion of induction. We proved the statement of theorem. □

When $h_1 = ... = h_n = h$, we assume

$$\partial^n f(x) \circ h^n = \partial^n f(x) \circ (h_1; ...; h_n)$$

This notation does not create ambiguity, because we can determine function according to number of arguments.

Theorem 3.2.4. *For any $n > 0$ following equation is true*

$$\partial^n p_n(x) \circ h^n = n! p_n(h)$$

Proof. For $n = 1$ following equation is true

$$\partial p_1(x) \circ h = \partial(a_0 x a_1) \circ h = a_0 h a_1 = 1! p_1(h)$$

Assume the statement is true for $n - 1$. Then according to theorem 3.2.1

$$\begin{aligned}\partial^n p_n(x) \circ h^n = (\partial^n p_{n-1}(x) \circ h^n) x a_n + (\partial^{n-1} p_{n-1}(x) \circ h^{n-1}) h a_n \\ + ... + (\partial^{n-1} p_{n-1}(x) \circ h^{n-1}) h a_n\end{aligned} \tag{3.2.3}$$

First term equal 0 according to theorem 3.2.2. The rest $n$ terms equal, and according to suggestion of induction from equation (3.2.3) it follows

$$\partial^n p_n(x) \circ h = n(\partial^{n-1} p_{n-1}(x) \circ h) h a_n = n(n-1)! p_{n-1}(h) h a_n = n! p_n(h)$$

Therefore, statement of theorem is true for any $n$. □

Let $p(x)$ be polynomial of power $n$.[3.2]

$$p(x) = p_0 + p_{1 i_1}(x) + ... + p_{n i_n}(x)$$

We assume sum by index $i_k$ which enumerates terms of power $k$. According to theorem 3.2.2, 3.2.3, 3.2.4

$$\partial^k p(0) \circ x = k! p_{k i_k}(x)$$

Therefore, we can write

$$p(x) = p_0 + (1!)^{-1} \partial p(0) \circ x + (2!)^{-1} \partial^2 p(0) \circ x^2 + ... + (n!)^{-1} \partial^n p(0) \circ x^n$$

This representation of polynomial is called **Taylor polynomial**. If we consider substitution of variable $x = y - y_0$, then considered above construction remain true for polynomial

$$p(y) = p_0 + p_{1 i_1}(y - y_0) + ... + p_{n i_n}(y - y_0)$$

Therefore

$$p(y) = p_0 + (1!)^{-1} \partial p(y_0) \circ (y - y_0) + (2!)^{-1} \partial^2 p(y_0) \circ (y - y_0)^2 + ... + (n!)^{-1} \partial^n p(y_0) \circ (y - y_0)^n$$

Assume that function $f(x)$ is differentiable in the Gâteaux sense at point $x_0$ up to any order.[3.3]

Theorem 3.2.5. *If function $f(x)$ holds*

$$f(x_0) = \partial f(x_0) \circ h = ... = \partial^n f(x_0) \circ h^n = 0$$

*then for $t \to 0$ expression $f(x + th)$ is infinitesimal of order higher than $n$ with respect to $t$*

$$f(x_0 + th) = o(t^n)$$

---

[3.2] I consider Taylor polynomial for polynomials by analogy with construction of Taylor polynomial in [5], p. 246.

[3.3] I explore construction of Taylor series by analogy with construction of Taylor series in [5], p. 248, 249.

Proof. When $n = 1$ this statement follows from equation (2.3.6).

Let statement be true for $n-1$. Map

$$f_1(x) = \partial f(x) \circ h$$

satisfies to condition

$$f_1(x_0) = \partial f_1(x_0) \circ h = ... = \partial^{n-1} f_1(x_0) \circ h^{n-1} = 0$$

According to suggestion of induction

$$f_1(x_0 + th) = o(t^{n-1})$$

Then equation (2.3.5) gets form

$$o(t^{n-1}) = \lim_{t \to 0,\ t \in R} (t^{-1} f(x + th))$$

Therefore,

$$f(x + th) = o(t^n)$$

□

Let us form polynomial

$$p(x) = f(x_0) + (1!)^{-1} \partial f(x_0) \circ (x - x_0) + ... + (n!)^{-1} \partial^n f(x_0) \circ (x - x_0)^n$$

According to theorem 3.2.5

$$f(x_0 + t(x - x_0)) - p(x_0 + t(x - x_0)) = o(t^n)$$

Therefore, polynomial $p(x)$ is good approximation of map $f(x)$.

If the mapping $f(x)$ has the Gâteaux derivative of any order, then passing to the limit $n \to \infty$, we get expansion into series

$$f(x) = \sum_{n=0}^{\infty} (n!)^{-1} \partial^n f(x_0) \circ (x - x_0)^n$$

which is called **Taylor series**.

CHAPTER 4

CHAPTER 5

# Index

CHAPTER 6

# Special Symbols and Notations

# Производная отображения банаховой алгебры

Александр Клейн

Aleks_Kleyn@MailAPS.org
http://AleksKleyn.dyndns-home.com:4080/
http://sites.google.com/site/AleksKleyn/
http://arxiv.org/a/kleyn_a_1
http://AleksKleyn.blogspot.com/

Аннотация. Пусть $A$ - банахова алгебра над коммутативным кольцом $D$. Отображение $f:A\to A$ дифференцируема по Гато, если

$$f(x+a)-f(x)=\partial f(x)\circ a+o(a)$$

где производная Гато $\partial f(x)$ отображения $f$ - линейное отображение приращения $a$ и $o$ - такое непрерывное отображение, что

$$\lim_{a\to 0}\frac{|o(a)|}{|a|}=0$$

Предполагая, что определена производная Гато $\partial^{n-1}f(x)$ порядка $n-1$, мы определим

$$\partial^n f(x)\circ(a_1\otimes...\otimes a_n)=\partial(\partial^{n-1}f(x)\circ(a_1\otimes...\otimes a_{n-1}))\circ a_n$$

производную Гато порядка $n$ отображения $f$. Если отображение $f(x)$ имеет все производные, то отображение $f(x)$ имеет разложение в ряд Тейлора

$$f(x)=\sum_{n=0}^{\infty}(n!)^{-1}\partial^n f(x_0)\circ(x-x_0)^n$$

# Оглавление

Глава 1

# Предисловие

## 1.1. Предисловие

В основе математического анализа лежит возможность линейного приближения к отображению, и основные построения математического анализа уходят корнями в линейную алгебру. Так как произведение в поле коммутативно, то линейная алгебра над полем относительно проста. При переходе к алгебре над коммутативным кольцом, некоторые утверждения линейной алгебры сохраняются, но появляются и новые утверждения, которые меняют ландшафт линейной алгебры.

Здесь я хочу обратить внимание на эволюцию, которую претерпело понятие производной со времён Ньютона. Когда мы изучаем функции одной переменной, то производная в заданной точке является числом

$$dx^2 = 2x\ dx$$

Когда мы изучаем функцию нескольких переменных, выясняется, что числа недостаточно. Производная становится вектором или градиентом

$$z = x^2 + y^3$$
$$dz = 2x\ dx + 3y^2\ dy$$

При изучении отображений векторных пространств мы впервые говорим о производной как об операторе

$$x = u\sin v \qquad y = u\cos v \qquad z = u$$
$$dx = \sin v\ du + u\cos v\ dv \qquad dy = \cos v\ du - u\sin v\ dv \qquad dz = 1\ du + 0\ dv$$

Но так как этот оператор линеен, то мы можем представить производную как матрицу. И в этом случае мы можем представить вектор приращения отображения как произведение матрицы производной (матрицы Якоби) на вектор приращения аргумента

$$\begin{pmatrix} dx \\ dy \\ dz \end{pmatrix} = \begin{pmatrix} \sin v & u\cos v \\ \cos v & -u\sin v \\ 1 & 0 \end{pmatrix} \begin{pmatrix} du \\ dv \end{pmatrix}$$

Поскольку алгебра является модулем над некоторым коммутативным кольцом, существует два пути изучения структур, порождённых над алгеброй.

Если алгебра является свободным модулем, то мы можем выбрать базис и рассматривать все операции в координатах относительно заданного базиса. Хотя базис может быть произвольным, мы можем выбрать наиболее простой базис с точки зрения алгебраических операций. Этот подход имеет без сомнения то преимущество, что мы работаем в коммутативном кольце, где все операции хорошо изучены.

Рассмотрение операции в алгебре независимо от выбранного базиса даёт возможность рассматривать элементы алгебры как самостоятельные объекты. Структура линейного отображения алгебры рассмотрена в книге [1]. В этой книге я использовал этот инструмент для изучения математического анализа над некоммутативной алгеброй.

## 1.2. Соглашения

Соглашение 1.2.1. *В выражении вида*

$$a_{i\cdot 0}xa_{i\cdot 1}$$

*предполагается сумма по индексу* $i$. □

Соглашение 1.2.2. *Если тензор* $a\in A^{\otimes(n+1)}$ *имеет разложение*

$$a=a_{i\cdot 0}\otimes a_{i\cdot 1}\otimes ...\otimes a_{i\cdot n}\quad i\in I$$

*то множество перестановок* $\sigma=\{\sigma_i\in S(n):i\in I\}$ *и тензор a порождают отображение*

$$(a,\sigma):A^{\times n}\to A$$

*определённое равенством*

$$\begin{aligned}(a,\sigma)\circ(b_1,...,b_n)&=(a_{i\cdot 0}\otimes a_{i\cdot 1}\otimes ...\otimes a_{i\cdot n},\sigma_i)\circ(b_1,...,b_n)\\&=a_{i\cdot 0}\sigma_i(b_1)a_{i\cdot 1}...\sigma_i(b_n)a_{i\cdot n}\end{aligned}$$

□

Соглашение 1.2.3. *Пусть тензор* $a\in A^{\otimes(n+1)}$ . *Если* $x_1=...=x_n=x$, *то мы положим*

$$a\circ x^n=a\circ(x_1\otimes ...\otimes x_n)$$

□

Соглашение 1.2.4. *Элемент* $\Omega$*-алгебры* $A$ *называется* $A$**-числом**. *Например, комплексное число также называется* $C$*-числом, а кватернион называется* $H$*-числом.* □

Соглашение 1.2.5. *Пусть* $A$ *ассоциативная* $D$*-алгебра. Представление*

$$A\otimes A\xrightarrow{f}{*} A\qquad f(p):a\to p\circ a$$

$D$*-модуля* $A\otimes A$ *определено равенством*

$$(a\otimes b)\circ c=acb \tag{1.2.1}$$

*и порождает множество линейных отображений. Это представление порождает произведение* $\circ$ *в* $D$*-модуле* $A\otimes A$ *согласно правилу*

$$(p\circ q)\circ a=p\circ(q\circ a)$$

□

Глава 2

# Дифференцируемые отображения

В этой главе мы изучим производную и дифференциал отображения в $D$ -алгебру. Поле комплексных чисел - это алгебра над полем действительных чисел. В математическом анализе функций комплексного переменного мы рассматриваем линейные отображения, порождённые отображением $I_0 \circ z = z$. Поэтому и в этом разделе мы ограничимся линейными отображениями, порождёнными отображением $I_0 \circ a = a$.

## 2.1. Топологическое кольцо

Определение 2.1.1. *Кольцо $D$ называется* **топологическим кольцом**[2.1]*, если $D$ является топологическим пространством, и алгебраические операции, определённые в $D$, непрерывны в топологическом пространстве $D$.* □

Согласно определению, для произвольных элементов $a, b \in D$ и для произвольных окрестностей $W_{a-b}$ элемента $a-b$, $W_{ab}$ элемента $ab$ существуют такие окрестности $W_a$ элемента $a$ и $W_b$ элемента $b$, что $W_a - W_b \subset W_{a-b}$, $W_a W_b \subset W_{ab}$.

Определение 2.1.2. **Норма на кольце** $D$ *- это отображение*[2.2]

$$d \in D \to |d| \in R$$

*такое, что*

- $|a| \geq 0$
- $|a| = 0$ *равносильно* $a = 0$
- $|ab| = |a|\ |b|$
- $|a + b| \leq |a| + |b|$

*Кольцо $D$, наделённое структурой, определяемой заданием на $D$ нормы, называется* **нормированным кольцом**. □

Инвариантное расстояние на аддитивной группе кольца $D$

$$d(a, b) = |a - b|$$

определяет топологию метрического пространства, согласующуюся со структурой кольца в $D$.

Определение 2.1.3. *Пусть $D$ - нормированное кольцо. Элемент $a \in D$ называется* **пределом последовательности** $\{a_n\}$

$$a =$$

[2.1] Определение дано согласно определению из [4], глава 4

[2.2] Определение дано согласно определению из [2], гл. IX, §3, п°2, а также согласно определению [6]-1.1.12, с. 23.

*если для любого* $\epsilon \in R$, $\epsilon > 0$, *существует, зависящее от* $\epsilon$, *натуральное число* $n_0$ *такое, что*

$$|a_n - a| < \epsilon$$

*для любого* $n > n_0$. □

Теорема 2.1.4. *Пусть* $D$ *- нормированное кольцо характеристики* 0 *и пусть* $d \in D$. *Пусть* $a \in D$ *- предел последовательности* $\{a_n\}$. *Тогда*

$$\lim_{n\to\infty} (a_n d) = ad$$

$$\lim_{n\to\infty} (d a_n) = da$$

Доказательство. Утверждение теоремы тривиально, однако я привожу доказательство для полноты текста. Поскольку $a \in D$ - предел последовательности $\{a_n\}$, то согласно определению 2.1.3 для заданного $\epsilon \in R$, $\epsilon > 0$, существует натуральное число $n_0$ такое, что

$$|a_n - a| < \frac{\epsilon}{|d|}$$

для любого $n > n_0$. Согласно определению 2.1.2 утверждение теоремы следует из неравенств

$$|a_n d - ad| = |(a_n - a)d| = |a_n - a||d| < \frac{\epsilon}{|d|}|d| = \epsilon$$

$$|d a_n - da| = |d(a_n - a)| = |d||a_n - a| < |d|\frac{\epsilon}{|d|} = \epsilon$$

для любого $n > n_0$. □

Определение 2.1.5. *Пусть* $D$ *- нормированное кольцо. Последовательность* $\{a_n\}$, $a_n \in D$ *называется* **фундаментальной** *или* **последовательностью Коши**, *если для любого* $\epsilon \in R$, $\epsilon > 0$ *существует, зависящее от* $\epsilon$, *натуральное число* $n_0$ *такое, что* $|a_p - a_q| < \epsilon$ *для любых* $p$, $q > n_0$. □

Определение 2.1.6. *Нормированное кольцо* $D$ *называется* **полным** *если любая фундаментальная последовательность элементов данного кольца сходится, т. е. имеет предел в этом кольце.* □

В дальнейшем, говоря о нормированном кольце характеристики 0, мы будем предполагать, что определён гомеоморфизм поля рациональных чисел $Q$ в кольцо $D$.

Теорема 2.1.7. *Полное кольцо* $D$ *характеристики* 0 *содержит в качестве подполя изоморфный образ поля* $R$ *действительных чисел. Это поле обычно отождествляют с* $R$.

Доказательство. Рассмотрим фундаментальную последовательность рациональных чисел $\{p_n\}$. Пусть $p'$ - предел этой последовательности в кольце $D$. Пусть $p$ - предел этой последовательности в поле $R$. Так как вложение поля $Q$ в тело $D$ гомеоморфно, то мы можем отождествить $p' \in D$ и $p \in R$. □

Теорема 2.1.8. *Пусть* $D$ *- полное кольцо характеристики* 0 *и пусть* $d \in D$. *Тогда любое действительное число* $p \in R$ *коммутирует с* $d$.

Доказательство. Мы можем представить действительное число $p \in R$ в виде фундаментальной последовательности рациональных чисел $\{p_n\}$. Утверждение теоремы следует из цепочки равенств

$$pd = \lim_{n\to\infty}(p_n d) = \lim_{n\to\infty}(dp_n) = dp$$

основанной на утверждении теоремы 2.1.4. □

## 2.2. Топологическая $D$-алгебра

Определение 2.2.1. *Пусть $D$ - топологическое коммутативное кольцо. $D$-алгебра $A$ называется* **топологической $D$-алгеброй**[2.3]*, если $A$ наделено топологией, согласующейся со структурой аддитивной группы в $A$, и отображения*

$$(a, v) \in D \times A \to av \in A$$
$$(v, w) \in A \times A \to vw \in A$$

*непрерывны.* □

Определение 2.2.2. **Норма на $D$-алгебре** *$A$ над нормированым коммутативным кольцом $D$*[2.4] *- это отображение*

$$a \in A \to |a| \in R$$

*такое, что*

- $|a| \geq 0$
- $|a| = 0$ *равносильно* $a = 0$
- $|a + b| \leq |a| + |b|$
- $|ab| = |a|\ |b|$
- $|da| = |d|\ |a|$, $d \in D$, $a \in A$

*$D$-алгебра $A$ над нормированным коммутативным кольцом $D$, наделённая структурой, определяемой заданием на $A$ нормы, называется* **нормированной $D$-алгеброй**. □

Определение 2.2.3. *Пусть $A$ - нормированная $D$-алгебра. Элемент $a \in A$ называется* **пределом последовательности** $\{a_n\}$

$$a = \lim_{n\to\infty} a_n$$

*если для любого $\epsilon \in R$, $\epsilon > 0$ существует, зависящее от $\epsilon$, натуральное число $n_0$ такое, что $|a_n - a| < \epsilon$ для любого $n > n_0$.* □

Определение 2.2.4. *Пусть $A$ - нормированная $D$-алгебра. Последовательность $\{a_n\}$, $a_n \in A$, называется* **фундаментальной** *или* **последовательностью Коши***, если для любого $\epsilon \in R$, $\epsilon > 0$, существует, зависящее от $\epsilon$, натуральное число $n_0$ такое, что $|a_p - a_q| < \epsilon$ для любых $p$, $q > n_0$.* □

Определение 2.2.5. *Нормированная $D$-алгебра $A$ называется* **банаховой $D$-алгеброй** *если любая фундаментальная последовательность элементов алгебры $A$ сходится, т. е. имеет предел в алгебре $A$.* □

Определение 2.2.6. *Пусть $A$ - банаховая $D$-алгебра. Множество элементов $a \in A$, $|a| = 1$, называется* **единичной сферой в алгебре** *$A$.* □

[2.3] Определение дано согласно определению из [3], с. 21

[2.4] Определение дано согласно определению из [2], гл. IX, §3, п°3

Определение 2.2.7. *Отображение*

$$f : A_1 \to A_2$$

*банаховой $D_1$-алгебры $A_1$ с нормой $|x|_1$ в банаховую $D_2$-алгебру $A_2$ с нормой $|y|_2$ называется* **непрерывным**, *если для любого сколь угодно малого $\epsilon > 0$ существует такое $\delta > 0$, что*

$$|x' - x|_1 < \delta$$

*влечёт*

$$|f(x') - f(x)|_2 < \epsilon$$

□

Определение 2.2.8. *Пусть*

$$f : A_1 \to A_2$$

*отображение банаховой $D_1$-алгебры $A_1$ с нормой $|x|_1$ в банаховую $D_2$-алгебру $A_2$ с нормой $|y|_2$. Величина*

$$\|f\| = sup\frac{|f(x)|_2}{|x|_1} \tag{2.2.1}$$

*называется* **нормой отображения** *$f$.* □

Теорема 2.2.9. *Пусть*

$$f : A_1 \to A_2$$

*линейное отображение банаховой $D_1$-алгебры $A_1$ с нормой $|x|_1$ в банаховую $D_2$-алгебру $A_2$ с нормой $|y|_2$. Тогда*

$$\|f\| = sup\{|f(x)|_2 : |x|_1 = 1\} \tag{2.2.2}$$

Доказательство. Из определений [1]-4.2.1, [1]-6.1.1 и теорем 2.1.7, 2.1.8 следует

$$f(rx) = rf(x) \quad r \in R \tag{2.2.3}$$

Из равенства (2.2.3) и определения 2.2.2 следует

$$\frac{|f(rx)|_2}{|rx|_1} = \frac{|r|\ |f(x)|_2}{|r|\ |x|_1} = \frac{|f(x)|_2}{|x|_1}$$

Полагая $r = \dfrac{1}{|x|_1}$, мы получим

$$\frac{|f(x)|_2}{|x|_1} = \left|f\left(\frac{x}{|x|_1}\right)\right|_2 \tag{2.2.4}$$

Равенство (2.2.2) следует из равенств (2.2.4) и (2.2.1). □

Теорема 2.2.10. *Пусть*

$$f : A_1 \to A_2$$

*линейное отображение банаховой $D_1$-алгебры $A_1$ с нормой $|x|_1$ в банаховую $D_2$-алгебру $A_2$ с нормой $|y|_2$. Отображение $f$ непрерывно, если $\|f\| < \infty$.*

Доказательство. Поскольку отображение $f$ линейно, то согласно определению 2.2.8

$$|f(x)-f(y)|_2=|f(x-y)|_2\le\|f\|\ |x-y|_1$$

Возьмём произвольное $\epsilon>0$. Положим $\delta=\frac{\epsilon}{\|f\|}$. Тогда из неравенства

$$|x-y|_1<\delta$$

следует

$$|f(x)-f(y)|_2\le\|f\|\ \delta=\epsilon$$

Согласно определению 2.2.7 отображение $f$ непрерывно. □

## 2.3. Производная отображений алгебры

Определение 2.3.1. *Пусть $A$ - банаховая $D$-алгебра. Отображение*

$$f:A\to A$$

**дифференцируемо по Гато** *на множестве $U\subset A$, если в каждой точке $x\in U$ изменение отображения $f$ может быть представлено в виде*

$$f(x+a)-f(x)=\partial f(x)\circ a+o(a)=\frac{\partial f(x)}{\partial x}\circ a+o(a) \tag{2.3.1}$$

*где* **производная Гато** $\partial f(x)$ **отображения** $f$ *- линейное отображение приращения $a$ и $o:A\to A$ - такое непрерывное отображение, что*

$$\lim_{a\to 0}\frac{|o(a)|}{|a|}=0$$

□

Замечание 2.3.2. *Согласно определению 2.3.1 при заданном $x$ производная Гато $\partial f(x)\in\mathcal{L}(A;A)$. Следовательно, производная Гато отображения $f$ является отображением*

$$\partial f:A\to\mathcal{L}(A;A)$$

*Выражения $\partial f(x)$ и $\frac{\partial f(x)}{\partial x}$ являются разными обозначениями одного и того же отображения. Мы будем пользоваться обозначением $\frac{\partial f(x)}{\partial x}$, если хотим подчеркнуть, что мы берём производную Гато по переменной $x$.* □

Теорема 2.3.3. *Мы можем представить производную Гато отображения $f$ в виде*[2.5]

$$\begin{aligned}\frac{\partial f(x)}{\partial x}\circ a&=\left(\frac{\partial_{s\cdot 0}f(x)}{\partial x}\otimes\frac{\partial_{s\cdot 1}f(x)}{\partial x}\right)\circ a\\&=\frac{\partial_{s\cdot 0}f(x)}{\partial x}a\frac{\partial_{s\cdot 1}f(x)}{\partial x}\end{aligned} \tag{2.3.3}$$

[2.5] Формально, мы должны записать производную отображения в виде

$$\frac{\partial f(x)}{\partial x}\circ a=\left(\frac{\partial_{k\cdot s\cdot 0}f(x)}{\partial x}\otimes\frac{\partial_{k\cdot s\cdot 1}f(x)}{\partial x}\right)\circ I_k\circ a=\frac{\partial_{k\cdot s\cdot 0}f(x)}{\partial x}(I_k\circ a)\frac{\partial_{k\cdot s\cdot 1}f(x)}{\partial x} \tag{2.3.2}$$

Однако, например, в теории функций комплексного переменного рассматриваются только линейные отображения, порождённые отображением $I_0\circ z=z$. Поэтому при изучении производных мы также ограничимся линейные отображениями, порождёнными отображением $I_0$. Переход к общему случаю не составляет особого труда.

*Выражение* $\dfrac{\partial_{s\cdot p} f(x)}{\partial x}$, $p=0,\ 1,$ *называется* **компонентой производной Гато** *отображения* $f(x)$.

Доказательство. Следствие определения 2.3.1 и теоремы [1]-6.4.5. □

Из определений [1]-4.2.1, [1]-6.1.1 , 2.3.1 и теоремы 2.1.7 следует

(2.3.4) $$\partial f(x)\circ(ra)=r\partial f(x)\circ a$$

$$r\in R\quad r\neq 0\quad a\in A\quad a\neq 0$$

Комбинируя равенство (2.3.4) и определение 2.3.1, мы получим знакомое определение производной Гато

(2.3.5) $$\partial f(x)\circ a=\lim_{t\to 0,\ t\in R}(t^{-1}(f(x+ta)-f(x)))$$

Определения производной Гато (2.3.1) и (2.3.5) эквивалентны. На основе этой эквивалентности мы будем говорить, что отображение $f$ дифференцируемо по Гато на множестве $U\subset D$, если в каждой точке $x\in U$ изменение отображения $f$ может быть представлено в виде

(2.3.6) $$f(x+ta)-f(x)=t\partial f(x)\circ a+o(t)$$

где $o:R\to A$ - такое непрерывное отображение, что

$$\lim_{t\to 0}\frac{|o(t)|}{|t|}=0$$

Если бесконечно малая $a$ в равенстве (2.3.1) является дифференциалом $dx$, то равенство (2.3.1) становится определением **дифференциала Гато**

(2.3.7) $$\begin{aligned}\frac{\partial f(x)}{\partial x}\circ dx&=\left(\frac{\partial_{s\cdot 0} f(x)}{\partial x}\otimes\frac{\partial_{s\cdot 1} f(x)}{\partial x}\right)\circ dx\\&=\frac{\partial_{s\cdot 0} f(x)}{\partial x}dx\frac{\partial_{s\cdot 1} f(x)}{\partial x}\end{aligned}$$

Теорема 2.3.4. *Пусть $A$ - банахова $D$-алгебра. Пусть $\overline{\overline{e}}$ - базис алгебры $A$ над кольцом $D$.* **Стандартное представление производной Гато отображения**

$$f:A\to A$$

*имеет вид*

(2.3.8) $$\frac{\partial f(x)}{\partial x}=\frac{\partial^{ij} f(x)}{\partial x}\overline{e}_i\otimes\overline{e}_j$$

*Выражение* $\dfrac{\partial^{ij} f(x)}{\partial x}$ *в равенстве* (2.3.8) *называется* **стандартной компонентой производной Гато** *отображения* $f$.

Доказательство. Утверждение теоремы является следствием утверждения [1]-6.4.5.2. □

Теорема 2.3.5. *Пусть $A$ - банахова $D$-алгебра. Пусть $\overline{\overline{e}}$ - алгебры $A$ над кольцом $D$. Тогда производная Гато отображения*

$$f:D\to D$$

*можно записать в виде*

$$\frac{\partial f(x)}{\partial x}\circ dx = dx^i\frac{\partial f^j}{\partial x^i}\overline{e}_j \tag{2.3.9}$$

*где $dx\in A$ имеет разложение*

$$dx = dx^i\ e_i \qquad dx^i\in D$$

*относительно базиса $\overline{\overline{e}}$ и матрица Якоби отображения $f$ имеет вид*

$$\frac{\partial f^j}{\partial x^i} = \frac{\partial^{kr} f(x)}{\partial x}C^p_{ki}C^j_{pr} \tag{2.3.10}$$

Доказательство. Утверждение теоремы является следствием теоремы [1]-6.4.5. □

Теорема 2.3.6. *Пусть $A$ - банаховая $D$-алгебра. Пусть $f$, $g$ - дифференцируемые отображения*

$$f:A\to A\quad g:A\to A$$

*Отображение*

$$f+g:A\to A$$

*дифференцируемо и производная Гато удовлетворяет соотношению*

$$\partial(f+g)(x)=\partial f(x)+\partial g(x) \tag{2.3.11}$$

Доказательство. Согласно определению (2.3.5),

$$\begin{aligned}\partial(f+g)(x)\circ a &= \lim_{t\to 0,\ t\in R}(t^{-1}((f+g)(x+ta)-(f+g)(x)))\\ &= \lim_{t\to 0,\ t\in R}(t^{-1}(f(x+ta)+g(x+ta)-f(x)-g(x)))\\ &= \lim_{t\to 0,\ t\in R}(t^{-1}(f(x+ta)-f(x)))\\ &+ \lim_{t\to 0,\ t\in R}(t^{-1}(g(x+ta)-g(x)))\\ &= \partial f(x)\circ a+\partial g(x)\circ a\end{aligned} \tag{2.3.12}$$

Равенство (2.3.11) следует из равенства (2.3.12). □

Теорема 2.3.7. *Пусть $A$ - банаховая $D$-алгебра. Пусть*

$$h:A\times A\to A$$

*непрерывное билинейное отображение. Пусть $f$, $g$ - дифференцируемые отображения*

$$f:A\to A\quad g:A\to A$$

*Отображение*

$$h(f,g):A\to A$$

*дифференцируемо и дифференциал Гато удовлетворяет соотношению*

$$\partial h(f(x),g(x))\circ dx = h(\partial f(x)\circ dx, g(x))+h(f(x),\partial g(x)\circ dx) \tag{2.3.13}$$

Доказательство. Равенство (2.3.13) следует из цепочки равенств

$$\begin{aligned}\partial h(f(x),g(x))\circ a &= \lim_{t\to 0}(t^{-1}(h(f(x+ta),g(x+ta))-h(f(x),g(x))))\\ &= \lim_{t\to 0}(t^{-1}(h(f(x+ta),g(x+ta))-h(f(x),g(x+ta))))\\ &+ \lim_{t\to 0}(t^{-1}(h(f(x),g(x+ta))-h(f(x),g(x))))\\ &= h(\lim_{t\to 0}t^{-1}(f(x+ta)-f(x)),g(x))\\ &+ h(f(x),\lim_{t\to 0}t^{-1}(g(x+ta)-g(x)))\end{aligned}$$

основанной на определении (2.3.5). □

Соглашение 2.3.8. *Для заданного билинейного отображения*

$$h:A\times A\to A$$

*мы рассмотрим отображения*

$$h_1:\mathcal{L}(D;A;A)\times A\to\mathcal{L}(D;A;A)$$

$$h_2:A\times\mathcal{L}(D;A;A)\to\mathcal{L}(D;A;A)$$

*определённые равенствами*

$$h_1(f,v)\circ u=h(f\circ u,v)$$

$$h_2(u,f)\circ v=h(u,f\circ v)$$

*Мы будем пользоваться буквой $h$ для обозначения отображений $h_1$, $h_2$.* □

Теорема 2.3.9. *Пусть $A$ - банаховая $D$-алгебра. Пусть*

$$h:A\times A\to A$$

*непрерывное билинейное отображение. Пусть $f$, $g$ - дифференцируемые отображения*

$$f:A\to A\quad g:A\to A$$

*Отображение*

$$h(f,g):A\to A$$

*дифференцируемо и производная Гато удовлетворяет соотношению*

$$\partial h(f(x),g(x))=h(\partial f(x),g(x))+h(f(x),\partial g(x))\tag{2.3.14}$$

Доказательство. Равенство (2.3.14) следует из равенства (2.3.13) и соглашения 2.3.8. □

Теорема 2.3.10. *Пусть $A$ - банаховая $D$-алгебра. Пусть $f$, $g$ - дифференцируемые отображения*

$$f:A\to A\quad g:A\to A$$

*Производная Гато удовлетворяет соотношению*

$$\frac{\partial f(x)g(x)}{\partial x}\circ dx=\left(\frac{\partial f(x)}{\partial x}\circ dx\right)\,g(x)+f(x)\,\left(\frac{\partial g(x)}{\partial x}\circ dx\right)\tag{2.3.15}$$

$$\frac{\partial f(x)g(x)}{\partial x}=\frac{\partial f(x)}{\partial x}\,g(x)+f(x)\,\frac{\partial g(x)}{\partial x}\tag{2.3.16}$$

Доказательство. Теорема является следствием теорем 2.3.7, 2.3.9 и определения [1]-5.1.1. □

Теорема 2.3.11. *Пусть $A$ - банаховая $D$-алгебра. Допустим производная Гато отображения*

$$f:A\to A$$

*имеет разложение*

$$\frac{\partial f(x)}{\partial x}=\frac{\partial_{s\cdot 0}f(x)}{\partial x}\otimes\frac{\partial_{s\cdot 1}f(x)}{\partial x}\tag{2.3.17}$$

*Допустим производная Гато отображения $g:D\to D$ имеет разложение*

$$\frac{\partial g(x)}{\partial x}=\frac{\partial_{t\cdot 0}g(x)}{\partial x}\otimes\frac{\partial_{t\cdot 1}g(x)}{\partial x}\tag{2.3.18}$$

*Производная Гато отображения $f(x)g(x)$ имеет вид*

$$\frac{\partial f(x)g(x)}{\partial x}=\frac{\partial_{s\cdot 0}f(x)}{\partial x}\otimes\left(\frac{\partial_{s\cdot 1}f(x)}{\partial x}g(x)\right)+\left(f(x)\frac{\partial_{t\cdot 0}g(x)}{\partial x}\right)\otimes\frac{\partial_{t\cdot 1}g(x)}{\partial x}\tag{2.3.19}$$

$$\begin{aligned}&\frac{\partial_{s\cdot 0}f(x)g(x)}{\partial x}=\frac{\partial_{s\cdot 0}f(x)}{\partial x} &&\frac{\partial_{t\cdot 0}f(x)g(x)}{\partial x}=f(x)\frac{\partial_{t\cdot 0}g(x)}{\partial x}\\ &\frac{\partial_{s\cdot 1}f(x)g(x)}{\partial x}=\frac{\partial_{s\cdot 1}f(x)}{\partial x}g(x) &&\frac{\partial_{t\cdot 1}f(x)g(x)}{\partial x}=\frac{\partial_{t\cdot 1}g(x)}{\partial x}\end{aligned}\tag{2.3.20}$$

Доказательство. Подставим (2.3.17) и (2.3.18) в равенство (2.3.16)

$$\begin{aligned}\frac{\partial f(x)g(x)}{\partial x}\circ a&=\left(\frac{\partial f(x)}{\partial x}\circ a\right)\ g(x)+f(x)\ \left(\frac{\partial g(x)}{\partial x}\circ a\right)\\ &=\frac{\partial_{s\cdot 0}f(x)}{\partial x}a\frac{\partial_{s\cdot 1}f(x)}{\partial x}g(x)+f(x)\frac{\partial_{t\cdot 0}g(x)}{\partial x}a\frac{\partial_{t\cdot 1}g(x)}{\partial x}\end{aligned}\tag{2.3.21}$$

Опираясь на (2.3.21), мы определяем равенства (2.3.20). □

Теорема 2.3.12. *Пусть $A$ - банаховая $D$-алгебра. Пусть $f$, $g$ - дифференцируемые отображения*

$$f:A\to A\quad g:A\to A$$

*Производная Гато удовлетворяет соотношению*

$$\frac{\partial f(x)\otimes g(x)}{\partial x}\circ dx=\left(\frac{\partial f(x)}{\partial x}\circ dx\right)\otimes g(x)+f(x)\otimes\left(\frac{\partial g(x)}{\partial x}\circ dx\right)\tag{2.3.22}$$

$$\frac{\partial f(x)\otimes g(x)}{\partial x}=\frac{\partial f(x)}{\partial x}\otimes g(x)+f(x)\otimes\frac{\partial g(x)}{\partial x}\tag{2.3.23}$$

Доказательство. Теорема является следствием теорем 2.3.7, 2.3.9, [1]-4.4.5 и определения [1]-5.1.1. □

Замечание 2.3.13. *Пусть*

$$\frac{\partial f(x)}{\partial x}=\frac{\partial_{s\cdot 0}f(x)}{\partial x}\otimes\frac{\partial_{s\cdot 1}f(x)}{\partial x}\tag{2.3.24}$$

$$\frac{\partial g(x)}{\partial x}=\frac{\partial_{t\cdot 0}g(x)}{\partial x}\otimes\frac{\partial_{t\cdot 1}g(x)}{\partial x}\tag{2.3.25}$$

*Тогда*

$$\frac{\partial f(x)\otimes g(x)}{\partial x}=\frac{\partial_{s\cdot 0}f(x)}{\partial x}\otimes\frac{\partial_{s\cdot 1}f(x)}{\partial x}\otimes g(x)+f(x)\otimes\frac{\partial_{t\cdot 0}g(x)}{\partial x}\otimes\frac{\partial_{t\cdot 1}g(x)}{\partial x} \tag{2.3.26}$$

*Мы не пишем скобки, так как тензорное произведение ассоциативно и дистрибутивно относительно сложения (теоремы [1]-3.3.11, [1]-3.4.5 ).* □

Теорема 2.3.14. *Пусть $A$ - банахова $D$-алгебра. Если производная Гато $\partial f(x)$ существует в точке $x$ и имеет конечную норму, то отображение $f$ непрерывно в точке $x$.*

Доказательство. Из определения 2.2.8 следует

$$|\partial f(x)\circ a|\le\|\partial f(x)\|\ |a| \tag{2.3.27}$$

Из (2.3.1), (2.3.27) следует

$$|f(x+a)-f(x)|<|a|\ \|\partial f(x)\| \tag{2.3.28}$$

Возьмём произвольное $\epsilon>0$. Положим

$$\delta=\frac{\epsilon}{\|\partial f(x)\|}$$

Тогда из неравенства

$$|a|<\delta$$

следует

$$|f(x+a)-f(x)|\le\|\partial f(x)\|\ \delta=\epsilon$$

Согласно определению 2.2.7 отображение $f$ непрерывно в точке $x$. □

Теорема 2.3.15. *Пусть $A$ - банахова $D$-алгебра. Пусть отображение*

$$f:A\to A$$

*дифференцируемо по Гато в точке $x$. Тогда*

$$\partial f(x)\circ 0=0$$

Доказательство. Следствие определения 2.3.1 и теоремы [1]-4.2.5. □

Теорема 2.3.16. *Пусть $A$ - банахова $D$-алгебра. Пусть отображение*

$$f:A\to A$$

*дифференцируемо по Гато в точке $x$ и норма производной Гато отображения $f$ конечна*

$$\|\partial f(x)\|=F\le\infty \tag{2.3.29}$$

*Пусть отображение*

$$g:A\to A$$

*дифференцируемо по Гато в точке*

$$y=f(x) \tag{2.3.30}$$

*и норма производной Гато отображения $g$ конечна*

$$\|\partial g(y)\|=G\le\infty \tag{2.3.31}$$

*Отображение*

$$(g\circ f)(x)=g(f(x))$$

*дифференцируемо по Гато в точке* $x$

(2.3.32)
$$\begin{cases}\partial(g\circ f)(x)=\partial g(y)\circ\partial f(x)\\ \partial(g\circ f)(x)\circ a=\partial g(y)\circ\partial f(x)\circ a\end{cases}$$

(2.3.33)
$$\begin{cases}\dfrac{\partial_{st\cdot 0}(g\circ f)(x)}{\partial x}=\dfrac{\partial_{s\cdot 0}g(f(x))}{\partial f(x)}\dfrac{\partial_{t\cdot 0}f(x)}{\partial x}\\ \dfrac{\partial_{st\cdot 1}(g\circ f)(x)}{\partial x}=\dfrac{\partial_{t\cdot 1}f(x)}{\partial x}\dfrac{\partial_{s\cdot 1}g(f(x))}{\partial f(x)}\end{cases}$$

Доказательство. Согласно определению 2.3.1

(2.3.34)
$$g(y+b)-g(y)=\partial g(y)\circ b+o_1(b)$$

где $o_1:A\to A$ - такое непрерывное отображение, что

$$\lim_{b\to 0}\frac{|o_1(b)|}{|b|}=0$$

Согласно определению 2.3.1

(2.3.35)
$$f(x+a)-f(x)=\partial f(x)\circ a+o_2(a)$$

где $o_2:A\to A$ - такое непрерывное отображение, что

$$\lim_{a\to 0}\frac{|o_2(a)|}{|a|}=0$$

Согласно (2.3.35) смещение $a$ значения $x\in A$ приводит к смещению

(2.3.36)
$$b=\partial f(x)\circ a+o_2(a)$$

значения $y$. Используя (2.3.30), (2.3.36) в равенстве (2.3.34), мы получим

(2.3.37)
$$\begin{aligned}&g(f(x+a))-g(f(x))\\ =&g(f(x)+\partial f(x)\circ a+o_2(a))-g(f(x))\\ =&\partial g(f(x))\circ(\partial f(x)\circ a+o_2(a))-o_1(\partial f(x)\circ a+o_2(a))\end{aligned}$$

Согласно определениям 2.3.1, [1]-4.2.1, [1]-6.1.1 , из равенства (2.3.37) следует

(2.3.38)
$$\begin{aligned}&g(f(x+a))-g(f(x))\\ =&\partial g(f(x))\circ\partial f(x)\circ a+\partial g(f(x))\circ o_2(a)-o_1(\partial f(x)\circ a+o_2(a))\end{aligned}$$

Согласно определению 2.2.2

(2.3.39)
$$\begin{aligned}&\lim_{a\to 0}\frac{|\partial g(f(x))\circ o_2(a)-o_1(\partial f(x)\circ a+o_2(a))|}{|a|}\\ \le&\lim_{a\to 0}\frac{|\partial g(f(x))\circ o_2(a)|}{|a|}+\lim_{a\to 0}\frac{|o_1(\partial f(x)\circ a+o_2(a))|}{|a|}\end{aligned}$$

Из (2.3.31) следует

(2.3.40)
$$\lim_{a\to 0}\frac{|\partial g(f(x))\circ o_2(a)|}{|a|}\le G\lim_{a\to 0}\frac{|o_2(a)|}{|a|}=0$$

Из (2.3.29) следует

$$\begin{aligned}&\lim_{a\to 0}\frac{|o_1(\partial f(x)\circ a+o_2(a))|}{|a|}\\&=\lim_{a\to 0}\frac{|o_1(\partial f(x)\circ a+o_2(a))|}{|\partial f(x)\circ a+o_2(a)|}\lim_{a\to 0}\frac{|\partial f(x)\circ a+o_2(a)|}{|a|}\\&\le\lim_{a\to 0}\frac{|o_1(\partial f(x)\circ a+o_2(a))|}{|\partial f(x)\circ a+o_2(a)|_2}\lim_{a\to 0}\frac{\|\partial f(x)\||a|+|o_2(a)|}{|a|}\\&=\lim_{a\to 0}\frac{|o_1(\partial f(x)\circ a+o_2(a))|}{|\partial f(x)\circ a+o_2(a)|}\|\partial f(x)\|\end{aligned}$$

Согласно теореме 2.3.15

$$\lim_{a\to 0}(\partial f(x)\circ a)+o_2(a))=0$$

Следовательно,

$$\lim_{a\to 0}\frac{|o_1(\partial f(x)\circ a+o_2(a))|}{|a|}=0 \tag{2.3.41}$$

Из равенств (2.3.39), (2.3.40), (2.3.41) следует

$$\lim_{a\to 0}\frac{|\partial g(f(x))\circ o_2(a)-o_1(\partial f(x)\circ a+o_2(a))|}{|a|}=0 \tag{2.3.42}$$

Согласно определению 2.3.1

$$(g\circ f)(x+a)-(g\circ f)(x)=\partial(g\circ f)(x)\circ a+o(a) \tag{2.3.43}$$

где $o:A\to A$ - такое непрерывное отображение, что

$$\lim_{a\to 0}\frac{|o(a)|}{|a|}=0$$

Равенство (2.3.32) следует из (2.3.38), (2.3.42), (2.3.43).

Из равенства (2.3.32) и теоремы 2.3.3 следует

$$\begin{aligned}&\frac{\partial_{st\cdot 0}(g\circ f)(x)}{\partial x}\;a\;\frac{\partial_{st\cdot 1}(g\circ f)(x)}{\partial x}\\&=\frac{\partial_{s\cdot 0}g(f(x))}{\partial f(x)}\;(\partial f(x)\circ a)\;\frac{\partial_{s\cdot 1}g(f(x))}{\partial f(x)}\\&=\frac{\partial_{s\cdot 0}g(f(x))}{\partial f(x)}\frac{\partial_{t\cdot 0}f(x)}{\partial x}\;a\;\frac{\partial_{t\cdot 1}f(x)}{\partial x}\;\frac{\partial_{s\cdot 1}g(f(x))}{\partial f(x)}\end{aligned} \tag{2.3.44}$$

(2.3.33) следуют из равенства (2.3.44). □

Глава 3

# Производная второго порядка отображения $D$-алгебры

## 3.1. Производная второго порядка отображения $D$-алгебры

Пусть $D$ - полное коммутативное кольцо характеристики 0. Пусть $A$ - ассоциативная $D$-алгебра. Пусть

$$f : A \to A$$

функция, дифференцируемая по Гато. Согласно замечанию 2.3.2 производная Гато является отображением

$$\partial f : A \to \mathcal{L}(A; A)$$

Согласно теореме [1]-6.2.5 и определению 2.2.8, множество $\mathcal{L}(A; A)$ является банаховой $D$-алгеброй. Следовательно, мы можем рассмотреть вопрос, является ли отображение $\partial f$ дифференцируемым по Гато.

Согласно определению 2.3.1

$$(\partial f \circ (x + a_2)) \circ a_1 - (\partial f \circ x) \circ a_1 = \partial(\partial f(x) \circ a_1) \circ a_2 + o_2(a_2) \tag{3.1.1}$$

где $o_2 : A \to \mathcal{L}(A; A)$ - такое непрерывное отображение, что

$$\lim_{a_2 \to 0} \frac{\|o_2(a_2)\|}{|a_2|} = 0$$

Согласно определению 2.3.1 отображение $\partial(\partial f(x) \circ a_1) \circ a_2$ линейно по переменной $a_2$. Из равенства (3.1.1) следует, что отображение $\partial(\partial f(x) \circ a_1) \circ a_2$ линейно по переменной $a_1$. Следовательно, отображение $\partial(\partial f(x) \circ a_1) \circ a_2$ билинейно.

Определение 3.1.1. *Полилинейное отображение*

$$\partial^2 f(x) \circ (a_1; a_2) = \frac{\partial^2 f(x)}{\partial x^2} \circ (a_1; a_2) = \partial(\partial f(x) \circ a_1) \circ a_2 \tag{3.1.2}$$

*называется* **производной Гато второго порядка** *отображения $f$.* □

Замечание 3.1.2. *Согласно определению 3.1.1 при заданном $x$ производная Гато второго порядка $\partial^2 f(x) \in \mathcal{L}(A, A; A)$. Следовательно, производная Гато второго порядка отображения $f$ является отображением*

$$\partial^2 f : A \to \mathcal{L}(A, A; A)$$

*Согласно теореме [1]-4.4.4, мы можем также рассматривать отображение*

$$\partial^2 f(x) \circ (a_1 \otimes a_2) = \partial^2 f(x) \circ (a_1; a_2)$$

*Тогда*

$$\partial^2 f(x) \in \mathcal{L}(A \otimes A; A)$$
$$\partial^2 f : A \to \mathcal{L}(A \otimes A; A)$$

*Мы будем пользоваться тем же символом для обозначения отображения, так как по характеру аргумента ясно о каком отображении идёт речь.* □

Теорема 3.1.3. *Мы можем представить производную Гато второго порядка отображения $f$ в виде*

$$\partial^2 f(x) \circ (a_1; a_2) = \left( \frac{\partial^2_{s\cdot 0} f(x)}{\partial x^2} \otimes \frac{\partial^2_{s\cdot 1} f(x)}{\partial x^2} \otimes \frac{\partial^2_{s\cdot 2} f(x)}{\partial x^2}, \sigma_s \right) \circ (a_1; a_2)$$
$$= \frac{\partial^2_{s\cdot 0} f(x)}{\partial x^2} \sigma_s(a_1) \frac{\partial^2_{s\cdot 1} f(x)}{\partial x^2} \sigma_s(a_2) \frac{\partial^2_{s\cdot 2} f(x)}{\partial x^2}$$

*Мы будем называть выражение* [3.1]

$$\frac{\partial^2_{s\cdot p} f(x)}{\partial x^2} \qquad p = 0, 1, 2$$

**компонентой производной Гато второго порядка** *отображения* $f(x)$.

Доказательство. Следствие определения 3.1.1 и теоремы [1]-6.6.6. □

По индукции, предполагая, что определена производная Гато $\partial^{n-1} f(x)$ порядка $n-1$, мы определим

(3.1.3)
$$\partial^n f(x) \circ (a_1; ...; a_n) = \frac{\partial^n f(x)}{\partial x^n} \circ (a_1; ...; a_n)$$
$$= \partial(\partial^{n-1} f(x) \circ (a_1; ...; a_{n-1})) \circ a_n$$

**производную Гато порядка** $n$ отображения $f$. Мы будем также полагать $\partial^0 f(x) = f(x)$.

## 3.2. Ряд Тейлора

Пусть $D$ - полное коммутативное кольцо характеристики 0. Пусть $A$ - ассоциативная $D$-алгебра. Пусть $p_k(x)$ - одночлен степени $k$, $k > 0$, одной переменной над $D$-алгеброй $A$.

Очевидно, что одночлен степени 0 имеет вид $a_0$, $a_0 \in A$. Для $k > 0$,

$$p_k(x) = p_{k-1}(x) x a_k$$

где $a_k \in A$. Действительно, последний множитель одночлена $p_k(x)$ является либо $a_k \in A$, либо имеет вид $x^l$, $l \ge 1$. В последнем случае мы положим $a_k = 1$. Множитель, предшествующий $a_k$, имеет вид $x^l$, $l \ge 1$. Мы можем представить этот множитель в виде $x^{l-1} x$. Следовательно, утверждение доказано.

В частности, одночлен степени 1 имеет вид $p_1(x) = a_0 x a_1$.

Не нарушая общности, мы можем положить $k = n$.

[3.1] Мы полагаем

$$\frac{\partial^2_{s\cdot p} f(x)}{\partial x^2} = \frac{\partial^2_{s\cdot p} f(x)}{\partial x \partial x}$$

Теорема 3.2.1. *Для произвольного $m>0$ справедливо равенство*

$$\begin{aligned}\partial^m(f(x)x)\circ(h_1;...;h_m)=&\partial^m f(x)\circ(h_1;...;h_m)x\\ &+\partial^{m-1}f(x)\circ(\widehat{h_1};...;h_{m-1};h_m)h_1+...\\ &+\partial^{m-1}f(x)\circ(h_1;...;h_{m-1};\widehat{h_m})h_m\end{aligned}\tag{3.2.1}$$

*где символ $\widehat{h^i}$ означает отсутствие переменной $h^i$ в списке.*

Доказательство. Для $m=1$ - это следствие равенства (2.3.16)

$$\partial(f(x)x)\circ h_1=(\partial f(x)\circ h_1)x+f(x)h_1$$

Допустим, (3.2.1) справедливо для $m-1$. Тогда

$$\begin{aligned}\partial^{m-1}(f(x)x)\circ(h_1;...;h_{m-1})=&\partial^{m-1}f(x)\circ(h_1;...;h_{m-1})x\\ &+\partial^{m-2}f(x)\circ(\widehat{h_1};...;h_{m-2};h_{m-1})h_1+...\\ &+\partial^{m-2}f(x)\circ(h_1;...;h_{m-2};\widehat{h_{m-1}})h_{m-1}\end{aligned}$$

Так как $\partial h_i=0$, то, пользуясь равенством (2.3.16), получим

$$\begin{aligned}\partial^m(f(x)x)\circ(h_1;...;h_{m-1};h_m)=&\partial^m f(x)\circ(h_1;...;h_{m-1};h_m)x\\ &+\partial^{m-1}f(x)\circ(h_1;...;h_{m-2};h_{m-1})h_m\\ &+\partial^{m-1}f(x)\circ(\widehat{h_1};...;h_{m-2};h_{m-1};h_m)h_1+...\\ &+\partial^{m-1}f(x)\circ(h_1;...;h_{m-2};\widehat{h_{m-1}};h_m)h_{m-1}\end{aligned}\tag{3.2.2}$$

Равенства (3.2.1) и (3.2.2) отличаются только формой записи. Теорема доказана. □

Теорема 3.2.2. *Для произвольного $n\ge 0$ справедливо равенство*

$$\partial^{n+1}p_n(x)=0$$

Доказательство. Так как $p_0(x)=a_0$, $a_0\in D$, то при $n=0$ теорема является следствием определения 2.3.1. Пусть утверждение теоремы верно для $n-1$. Согласно теореме 3.2.1, при условии $f(x)=p_{n-1}(x)$ мы имеем

$$\begin{aligned}\partial^{n+1}p_n(x)(h_1;...;h_{n+1})=&\partial^{n+1}(p_{n-1}(x)xa_n)(h_1;...;h_{n+1})\\ =&\partial^{n+1}p_{n-1}(x)(h_1;...;h_{n+1})xa_n\\ &+\partial^n p_{n-1}(x)(\widehat{h_1};...;h_n;h_{n+1})h_1a_n+...\\ &+\partial^n p_{n-1}(x)(h_1;...;h_n;\widehat{h_{n+1}})h_{n+1}a_n\end{aligned}$$

Согласно предположению индукции все одночлены равны 0. □

Теорема 3.2.3. *Если $m<n$, то справедливо равенство*

$$\partial^m p_n(0)=0$$

Доказательство. Для $n=1$ справедливо равенство

$$\partial^0 p_1(0)=a_0xa_1=0$$

Допустим, утверждение справедливо для $n-1$. Тогда согласно теореме 3.2.1

$$\begin{aligned}\partial^m(p_{n-1}(x)xa_n)(h_1;...;h_m)=&\partial^m p_{n-1}(x)(h_1;...;h_m)xa_n\\&+\partial^{m-1}p_{n-1}(x)(\widehat{h_1};...;h_{m-1};h_m)h_1a_n+...\\&+\partial^{m-1}p_{n-1}(x)(h_1;...;h_{m-1};\widehat{h_m})h_ma_n\end{aligned}$$

Первое слагаемое равно 0 так как $x=0$. Так как $m-1<n-1$, то остальные слагаемые равны 0 согласно предположению индукции. Утверждение теоремы доказано. □

Если $h_1=...=h_n=h$, то мы положим

$$\partial^n f(x)\circ h^n=\partial^n f(x)\circ(h_1;...;h_n)$$

Эта запись не будет приводить к неоднозначности, так как по числу аргументов ясно, о какой функции идёт речь.

Теорема 3.2.4. *Для произвольного $n>0$ справедливо равенство*

$$\partial^n p_n(x)\circ h^n=n!p_n(h)$$

Доказательство. Для $n=1$ справедливо равенство

$$\partial p_1(x)\circ h=\partial(a_0xa_1)\circ h=a_0ha_1=1!p_1(h)$$

Допустим, утверждение справедливо для $n-1$. Тогда согласно теореме 3.2.1

$$\begin{aligned}\partial^n p_n(x)\circ h^n=&(\partial^n p_{n-1}(x)\circ h^n)xa_n+(\partial^{n-1}p_{n-1}(x)\circ h^{n-1})ha_n\\&+...+(\partial^{n-1}p_{n-1}(x)\circ h^{n-1})ha_n\end{aligned}\tag{3.2.3}$$

Первое слагаемое равно 0 согласно теореме 3.2.2. Остальные $n$ слагаемых равны, и согласно предположению индукции из равенства (3.2.3) следует

$$\partial^n p_n(x)\circ h=n(\partial^{n-1}p_{n-1}(x)\circ h)ha_n=n(n-1)!p_{n-1}(h)ha_n=n!p_n(h)$$

Следовательно, утверждение теоремы верно для любого $n$. □

Пусть $p(x)$ - многочлен степени $n$.[3.2]

$$p(x)=p_0+p_{1i_1}(x)+...+p_{ni_n}(x)$$

Мы предполагаем сумму по индексу $i_k$, который нумерует слагаемые степени $k$. Согласно теоремам 3.2.2, 3.2.3, 3.2.4

$$\partial^k p(0)\circ x=k!p_{ki_k}(x)$$

Следовательно, мы можем записать

$$p(x)=p_0+(1!)^{-1}\partial p(0)\circ x+(2!)^{-1}\partial^2 p(0)\circ x^2+...+(n!)^{-1}\partial^n p(0)\circ x^n$$

Это представление многочлена называется **формулой Тейлора для многочлена**. Если рассмотреть замену переменных $x=y-y_0$, то рассмотренное построение остаётся верным для многочлена

$$p(y)=p_0+p_{1i_1}(y-y_0)+...+p_{ni_n}(y-y_0)$$

откуда следует

$$p(y)=p_0+(1!)^{-1}\partial p(y_0)\circ(y-y_0)+(2!)^{-1}\partial^2 p(y_0)\circ(y-y_0)^2+...+(n!)^{-1}\partial^n p(y_0)\circ(y-y_0)^n$$

[3.2] Я рассматриваю формулу Тейлора для многочлена по аналогии с построением формулы Тейлора в [5], с. 246.

Предположим, что функция $f(x)$ в точке $x_0$ дифференцируема в смысле Гато до любого порядка.[3.3]

Теорема 3.2.5. *Если для функции $f(x)$ выполняется условие*

$$f(x_0) = \partial f(x_0) \circ h = ... = \partial^n f(x_0) \circ h^n = 0$$

*то при $t \to 0$ выражение $f(x+th)$ является бесконечно малой порядка выше $n$ по сравнению с $t$*

$$f(x_0 + th) = o(t^n)$$

Доказательство. При $n = 1$ это утверждение следует из равенства (2.3.6).

Пусть утверждение справедливо для $n - 1$. Для отображения

$$f_1(x) = \partial f(x) \circ h$$

выполняется условие

$$f_1(x_0) = \partial f_1(x_0) \circ h = ... = \partial^{n-1} f_1(x_0) \circ h^{n-1} = 0$$

Согласно предположению индукции

$$f_1(x_0 + th) = o(t^{n-1})$$

Тогда равенство (2.3.5) примет вид

$$o(t^{n-1}) = \lim_{t \to 0,\ t \in R} (t^{-1} f(x + th))$$

Следовательно,

$$f(x + th) = o(t^n)$$

□

Составим многочлен

$$p(x) = f(x_0) + (1!)^{-1} \partial f(x_0) \circ (x - x_0) + ... + (n!)^{-1} \partial^n f(x_0) \circ (x - x_0)^n$$

Согласно теореме 3.2.5

$$f(x_0 + t(x - x_0)) - p(x_0 + t(x - x_0)) = o(t^n)$$

Следовательно, полином $p(x)$ является хорошей апроксимацией отображения $f(x)$.

Если отображение $f(x)$ имеет производную Гато любого порядка, то переходя к пределу $n \to \infty$, мы получим разложение в ряд

$$f(x) = \sum_{n=0}^{\infty} (n!)^{-1} \partial^n f(x_0) \circ (x - x_0)^n$$

который называется **рядом Тейлора**.

[3.3] Я рассматриваю построение ряда Тейлора по аналогии с построением ряда Тейлора в [5], с. 248, 249.

Глава 4

# Список литературы

Глава 5

# Предметный указатель

Глава 6

# Специальные символы и обозначения